\documentclass[aos,MSNbibl,seceqn,nameyear,dvips]{arximspdf}

\doi{10.1214/11-AOS923}
\volume{39}
\issue{6}
\pubyear{2011}
\firstpage{2852}
\lastpage{2882}

\makeatletter

\newtheorem{lemma}{Lemma}[section]
\newtheorem{corollary}{Corollary}[section]
\newtheorem{theorem}{Theorem}[section]

\newproclaim{remark}{Remark}[section]

\newcommand{\ba}{{\mathbf{a}}}
\newcommand{\bx}{{\mathbf{x}}}
\newcommand{\bU}{{\mathbf{U}}}

\newcommand{\bV}{{\mathbf{V}}}
\newcommand{\bY}{{\mathbf{Y}}}
\newcommand{\bZ}{{\mathbf{Z}}}
\newcommand{\wba}{\widehat{\ba}}
\newcommand{\bmu}{{\bolds\mu}}
\newcommand{\bGa}{{\bolds\Gamma}}
\newcommand{\bSi}{{\bolds\Sigma}}
\newcommand{\wpi}{\widehat{\pi}}

\newcommand{\wtlam}{\widetilde{\lambda}}
\newcommand{\wlam}{\widehat{\lambda}}
\newcommand{\wphi}{\widehat{\phi}}
\newcommand{\wtphi}{\widetilde{\phi}}
\newcommand{\trasp}{^{\mathrm{T}}}
\newcommand{\var}{\operatorname{var}}
\newcommand{\cov}{\operatorname{cov}}

\newcommand{\convpp}{\stackrel{\mathrm{a.s.}}{\longrightarrow}}
\newcommand{\convppp}{\stackrel{\mathit{a.s.}}{\longrightarrow}}
\newcommand{\convweak}{\stackrel{\omega}{\longrightarrow}}

\newcommand{\median}{\operatorname{median}}
\newcommand{\argmax}{\mathop{\arg\max}}
\newcommand{\argmin}{\mathop{\arg\min}}

\newcommand{\rob}{\mathrm{R}}

\newcommand{\raw}{\mathrm{RAW}}

\newcommand{\smooth}{{\mathrm{S}}}

\newcommand{\smooths}{{\mathrm{PS}}}
\newcommand{\smoothn}{{\mathrm{PN}}}
\newcommand{\sieve}{{\mathrm{SI}}}
\newcommand{\pr}{{\mathrm{PR}}}

\newcommand{\KCV}{{\mathrm{KCV}}}

\newcommand{\real}{\mathbb{R}}
\newcommand{\esp}{\mathbb{E}}
\newcommand{\prob}{\mathbb{P}}

\makeatother

\begin{document}
\begin{frontmatter}

\title{Robust functional principal components: A~projection-pursuit
approach}
\runtitle{Robust functional principal components}

\begin{aug}
\author[A]{\fnms{Juan Lucas} \snm{Bali}\thanksref{T1}\ead[label=e1]{lbali@dm.uba.ar}},
\author[A]{\fnms{Graciela} \snm{Boente}\corref{}\thanksref{T1,T4}\ead[label=e2]{gboente@dm.uba.ar}},\\
\author[B]{\fnms{David E.} \snm{Tyler}\thanksref{T2}\ead[label=e3]{dtyler@rci.rutgers.edu}}
\and
\author[C]{\fnms{Jane-Ling} \snm{Wang}\thanksref{T3,T4}\ead[label=e4]{wang@wald.ucdavis.edu}}
\runauthor{Bali, Boente, Tyler and Wang}
\affiliation{Universidad de Buenos Aires and CONICET,
Universidad de Buenos Aires and CONICET,
Rutgers University and University of California at Davis}
\address[A]{J. L. Bali\\
G. Boente\\
Instituto de C\'alculo\\
Facultad de Ciencias Exactas y Naturales\\
Ciudad Universitaria, Pabell\'on 2\\
Buenos Aires, 1428\\
Argentina\\
\printead{e1}\\
\phantom{E-mail: }\printead*{e2}}
\address[B]{D. E. Tyler\\
Department of Statistics \\
Hill Center, Busch Campus \\
Rutgers University \\
Piscataway New Jersey 08854 \\
USA\\
\printead{e3}}
\address[C]{J.-L. Wang\\
Department of Statistics\\
University of California \\
Davis, California 95616\\
USA \\
\printead{e4}} %adresu isvedimo komanda gale!
\end{aug}

\thankstext{T1}{Supported in part by Grants X018 from
Universidad of Buenos Aires, PID 112-200801-00216 from
CONICET and PICT 821 from ANPCYT,
Argentina.}

\thankstext{T2}{Supported in part by NSF Grant DMS-09-06773.}

\thankstext{T3}{Supported in part by NSF
Grant DMS-09-06813.}

\thankstext{T4}{This project was completed while
Graciela Boente and Jane-Ling Wang were visiting SAMSI in
2010. They would like to acknowledge the inspiring environment of the
2010--2011 AOOD SAMSI program as well as the support given from SAMSI.}

% HISTORY:
\received{\smonth{11} \syear{2010}}
\revised{\smonth{9} \syear{2011}}

% ABSTRACT
%
\begin{abstract}
In many situations, data are recorded over a period of time and may be
regarded as realizations of a stochastic process. In this paper, robust
estimators for the principal components are considered {by} adapting
the projection pursuit approach to the functional data setting. Our
approach combines robust projection-pursuit with different smoothing
methods. Consistency of the estimators are shown under mild
assumptions. The performance of the classical and robust procedures
are compared in a simulation study under different contamination
schemes.
\end{abstract}

% KEYWORDS
%

\end{frontmatter}

%s1 #&#
\section{Introduction}\label{intro}

Analogous to classical principal components analysis (PCA), the
projection-pursuit approach to robust PCA is based on finding projections
of the data which have maximal dispersion. Instead of using the
variance as a measure of dispersion, a robust scale estimator $s_{n}$
is used for the maximization problem. This approach was introduced by
\citet{li}, who proposed estimators based on maximizing (or
minimizing) a~robust scale. In this way, given a sample $\bx_i\in
\real
^d$, $1\le i\le n$, the first robust principal component vector is
defined as
\[
\wba= \argmax_{\{\ba\in\real^d\dvtx \ba\trasp\ba=1\}} s_{n}(\ba
\trasp\bx
_{1},\ldots,\ba\trasp\bx_{n}) .
\]
The subsequent principal component vectors are obtained by imposing
orthogonality conditions. In the multivariate setting, \citet{li}
argue that the breakdown point for this projection-pursuit based
procedure is the same as that of the scale estimator $s_{n}$. Later on,
\citet{crouruiz2} derived the influence functions of the resulting
principal components, while their asymptotic distribution was studied
in \citet{cui}. A maximization algorithm for obtaining $\wba$
was proposed in \citet{crouruiz} and adapted for
high-dimensional data in \citet{crfil}.

The aim of this paper is to adapt the projection pursuit approach to
the functional data setting. We focus on functional data that are
recorded over a period of time and regarded as realizations
of a stochastic process, often assumed to be in the $L^2$ space on a
real interval. Various choices of robust scales, including the median
of the absolute deviation {about the median} (\textsc{mad}) and
$M$-estimates of scale are considered and compared.

Classical functional PCA uses the eigenvalues and eigenfunctions of the
sample covariance operator. \citet{dpr} have studied the
asymptotic properties of these sample functional principal components.
\citet{ris} proposed to smooth the principal components
by imposing an additive roughness penalty to the sample variance. The
consistency of this method was subsequently studied by \citet{pez}. Another approach to smoothing the principal
components, proposed in \citet{s} and reviewed in \citet{rasi2},
is based on penalizing the norm rather than the
sample variance, while \citet{bo} considered a
kernel-based approach. More recent work on estimation of the principal
components and the covariance function includes \citet{ger}, \citet{hall},
\citet{hallmulwa} and
\citet{yao}.

The literature on robust principal components in the functional data
setting, though, is rather sparse. To our knowledge, the first attempt
to provide estimators of the principal components that are less
sensitive to anomalous observations was due to \citet{loc},
although their approach is multivariate in nature. \citet{ger08}
studied a fully functional approach to robust estimation of the
principal components by considering a functional version of the
spherical principal components defined in \citet{loc}.
\citet{hyn} give an application of a
robust projection-pursuit approach, applied to smoothed trajectories,
but did not study the properties of their method in detail.

In this paper, we introduce several robust estimators of the principal
components in the functional data setting. Our approach uses a robust
projection-pursuit combined with various smoothing methods. A primary
focus of this paper is to provide a rigorous theoretical foundation for
this approach to robust functional PCA. In particular, we establish
under very general conditions the strong consistency of the our
proposed estimators.

In Section \ref{prop}, the robust estimators of the principal
components, based on both the raw and smoothed approaches, are
introduced. Consistency results and the asymptotic robustness of the
procedure are established in Section~\ref{consist}, while
Fisher-consistency of the related functionals is studied in Section~\ref {fisher}.
Section~\ref{appen} provides conditions under which one
of the crucial assumptions hold. Selection of the smoothing parameters
for the smooth principal components is discussed in Section~\ref{smoothpar}.
The results of a Monte Carlo study are reported in Section~\ref{monte}.
Finally, Section~\ref{concl} contains some concluding
remarks. Most proofs are relegated to the \hyperref[appenA]{Appendix}
and to the technical supplementary material available online; see
\citet{bali3a}. We begin the next section with notation and a review of
some basic concepts which are utilized in later sections.

%s2 #&#
\section{Preliminaries}\label{prelim}
%s2.1 #&#
\subsection{Functional principal components analysis}\label{fpca}
$\!\!\!$Principal components ana\-lysis, which was originally developed for
multivariate data, has been successfully extended to accommodate
functional data, and is
usually referred to as functional PCA. Principal components analysis
for general Hilbert spaces can be described as follows.

Let $X\in\mathcal{H}$ be a random element of a Hilbert space
$\mathcal
H$ defined in $(\Omega,\mathcal{A},P)$. Denote by $\langle\cdot
,\cdot
\rangle$ the inner product in $\mathcal H$ and by $\|\alpha\|
^2=\langle
\alpha,\alpha\rangle$. Assume that~$X$ has finite second moment, that
is, $\esp(\|X\|^2)<\infty$. The bilinear operator $a_X\dvtx\mathcal{H}
\times\mathcal{H} \to\real$ defined as $a_X(\alpha,\beta) =
\cov(\langle\alpha,X\rangle, \langle\beta,X\rangle)$ leads to a
continuous operator. The Riesz representation theorem then implies that
there exists a bounded operator, $\bGa_X\dvtx\mathcal{H}\to\mathcal{H}$,
such that $a_X(\alpha,\beta) = \langle\alpha, \bGa_X \beta\rangle$.
The operator $\bGa_X$ is called the covariance operator of $X$ and is
linear, self-adjoint and continuous.

Although the general situation in which $X\in\mathcal{H}$ is treated
in this paper, to help simplify the
basic ideas, we first consider the case $X\in L^2(\mathcal{I})$ where
$\mathcal{I}\subset\real$ is a finite interval. We take the usual
inner product for
$L^2(\mathcal{I})$, that is, $\langle\alpha, \beta\rangle= \int
_\mathcal{I} \alpha(t) \beta(t)\,dt$ and denote the covariance function
of $X$ by
$\gamma_X(t,s) = \cov(X(t),X(s))$. The corresponding
covariance operator $\bGa_X\dvtx L^2(\mathcal{I})\to L^2(\mathcal{I})$ is
such that $\bGa_X(\alpha)(t)=\int_\mathcal{I} \gamma_X(t,s) \alpha(s)
\,ds$. It is assumed the covariance function satisfies
$\int_\mathcal{I} \int_\mathcal{I} \gamma_X^2(t,s) \,dt \,ds
<\infty
$. Consequently, $\bGa_X$ is a Hilbert--Schmidt operator.

A Hilbert--Schmidt operator has a countable number of eigenvalues, all
of which are real. $\mathcal{F}$ will stand for the Hilbert space of
such operators with inner product defined by $\langle\bGa_1, \bGa
_2\rangle_{\mathcal{F}} =\sum_{j=1}^\infty\langle\bGa_1 u_j, \bGa_2
u_j\rangle$, where $\{u_j \dvtx j\ge1\}$ is any orthonormal basis of
$L^2(\mathcal{I})$. Furthermore, since the covariance operator~$\bGa_X$
is also positive semi-definite, its eigenvalues are nonnegative. As
with symmetric matrices in finite-dimensional Euclidean spaces, one can
choose the eigenfunctions of a Hilbert--Schmidt operator so that they
form an orthonormal basis for $L^2(\mathcal{I})$. Let $\{\phi_j\dvtx j\ge
1\}$ and
$\{\lambda_j\dvtx j\ge1\}$ be respectively an orthonormal basis of
eigenfunctions and their corresponding eigenvalues for the covariance
operator~$\bGa_X$,
with $\lambda_j\ge\lambda_{j+1}$. The spectral value decomposition for
$\bGa_X$ can then be expressed as $\bGa_X = \sum_{j=1}^\infty
\lambda_j
\phi_j \otimes\phi_j$,
with $\otimes$ being the tensor product, or equivalently $\gamma_X(t,s)
= \sum_{j=1}^\infty\lambda_j \phi_j(t)\phi_j(s)$, with\vspace*{1pt}
$ \sum_{j=1}^\infty\lambda^2_j = \|\bGa_X\|_\mathcal{F}^2 = \int
_\mathcal{I}\int_\mathcal{I} \gamma_X^2(t,s)\,dt \,ds$. The $j$th principal
component variable is then defined as
$Z_j = \langle\phi_j, X \rangle$, which leads to the Karhunen--Lo\`
{e}ve expansion $X(t) =\mu(t)+ \sum_{j = 1}^{\infty} Z_j \phi_j(t)$,
with $\mu(t)=\esp(X(t))$ and the $Z_j$'s being uncorrelated and having
variances $\lambda_j$ in descending order.

In general, for $Y = \langle\alpha,X\rangle$, which is a linear
function of the process $\{X(s)\}$, we have $\var(Y) = \langle\alpha
,\bGa_X \alpha\rangle$. An important optimality
property of the first principal component variable is that it can be
defined as the variable $Z_1 = \langle\alpha_1,X\rangle$ such that
%
%e2.1 #&#
%
\begin{equation}\label{MAX}
\var(Z_1)=\sup_{\{\alpha\dvtx \|\alpha\|=1\}}
\var(\langle\alpha,X\rangle)=\sup_{\{\alpha\dvtx \|
\alpha\|=1\}}
\langle\alpha,\bGa_X \alpha\rangle.
\end{equation}
Any solution to (\ref{MAX}), that is, any $\alpha$ for which the
supremum is obtained, corresponds
to an eigenfunction associated with the largest eigenvalue of the
covariance operator $\bGa_X$, that is, $\alpha_1=\phi_1$ and $\var
(Z_1)=\lambda_1$. If
$\lambda_1 > \lambda_2$, then~$\alpha_1$ is unique up to a sign change.
As in the multivariate setting, the other principal components can be
obtained successively
via (\ref{MAX}), but under the orthogonality condition that $\langle
\alpha_j, \alpha_k \rangle=0$ for $j < k$.

%s2.2 #&#
\subsection{Scale functionals and estimates} \label{scale}
The basic idea underlying our approach is to view principal components
as in (\ref{MAX}), but with the variance replaced by a robust scale functional.
We first recall the definition of a~scale functional. Denote by
$\mathcal{G}$ the set of all univariate distributions. A scale
functional $\sigma_{\rob}\dvtx\mathcal{G} \to[0,+\infty)$
is one which is location invariant and scale equivariant, that is, if
$G_{a,b}$ stands for the distribution of $aY+b$ when $Y\sim G$, then,
$\sigma_{\rob}(G_{a,b})=|a|\sigma_{\rob}(G)$,
for all real numbers $a$ and $b$. Two well-known examples of scale
functionals are the standard deviation, $\mbox{\textsc{sd}}(G)=\{\esp
(Y-\esp(Y))^2\}^{1/2}$, where $Y\sim G$,
and the median absolute deviation about the median, $\mbox{\textsc{mad}}(G)= c
\median(|Y-\median(Y)| )$. The normalization constant $c$,
used in the \textsc{mad}, can be chosen so that its
empirical or sample version is consistent for a scale parameter of
interest. Typically, one chooses $c=1/\Phi^{-1}(0.75)$ so that the
\textsc{mad} equals the standard deviation at a normal distribution.

The breakdown points, a measure of global robustness, for the standard
deviation and the \textsc{mad} are $0$ and $1/2$, respectively. The
\textsc{mad}, however, has a discontinuous influence function, which
reflects some local instability. Furthermore, the empirical version of
the \textsc{mad} is known to be fairly inefficient at the normal and
other distributions; see \citet{hub}.
In the finite-dimensional setting, as reported in \citet{cui}
the impact of a discontinuous influence function on the efficiency of
the estimators of the principal directions is even more dramatic
covariance function.

A broader class of robust scale functionals, which includes as special
cases both the \textsc{sd} and the \textsc{mad}, are the $M$-scale
functionals. An $M$-scale
functional with a bounded and continuous score function can have both a
high breakdown point and a continuous and bounded influence function.
Also, their empirical versions,
the $M$-estimates of scale, can be tuned to have good efficiency over a
broad range of distributions. Given a location parameter $\mu$, an
$M$-scale functionals $\sigma_M(G)$
with a continuous score function $\chi\dvtx\real\to\real$ can be defined
to be a solution to the equation
%
%e2.2 #&#
%
\begin{equation}\label{sigmarob}
\esp\biggl[ \chi\biggl( \frac{ Y-\mu}{\sigma_{\rob}(G)} \biggr) \biggr]
=\delta.
\end{equation}
Given a location statistic $\widehat{\mu}_n$, the corresponding
$M$-estimate of scale is then a solution $\widehat{\sigma}_n$ to the
$M$-estimating equation
%
%e2.3 #&#
%
\begin{equation}\label{snrob}
\frac{1}{n}\sum_{i=1}^{n}\chi\biggl( \frac{Y_{i} -\widehat{\mu
}_n}{\widehat{\sigma}_{ n} } \biggr) = \delta.
\end{equation}
If the score function is discontinuous, as is the case with the \textsc
{mad}, then a slight modification to (\ref{sigmarob}) and (\ref{snrob})
is needed; see \citet{mart}.

Typically, the score function $\chi$ is even with $\chi(0) = 0$,
nondecreasing on $\real_+$ and with $0< \sup_{x\in\real}\chi
(x)=\chi
(+\infty)=\lim_{x\to+\infty} \chi(x)$.
When $\chi(+\infty)=2\delta$, the $M$-estimate of scale has a 50\%
breakdown point, and by choosing $\chi$ properly one can also obtain a
highly efficient estimate; see \citet{croux}.
One such popular choice, and the one we use in our Monte Carlo study,
is the score function introduced by \citet{beatuk}, namely
%
%e2.4 #&#
%
\begin{equation}\label{funcionBT}
\chi_c( y) = \min\bigl(3 ( y/c )^2 - 3 ( y/c)^4
+ ( y/c )^6, 1\bigr)
\end{equation}
with $c$ being a tuning constant chosen so that the corresponding
$M$-estimator of scale is consistent for a scale parameter of interest.
For example, the choice $c = 1.56$ when $\delta=1/2$ ensures
that the $M$-scale functional is Fisher-consistent at the normal
distribution and has a 50\% breakdown point.

For continuous and nondecreasing score functions $\chi$, the solutions
to (\ref{sigmarob}) and (\ref{snrob}) are unique, and the simple
re-weighting algorithm
\[
\bigl\{\widehat{\sigma}_{ n}^{(k+1)}\bigr\}^2=\frac{1}{n \delta}\sum
_{i=1}^{n} w
\biggl( \frac{Y_{i} -\widehat{\mu}}{\widehat{\sigma}_{n}^{(k)} }
\biggr)
(Y_{i} -\widehat{\mu})^2,
\]
where $w(y)=\chi(y)/y^2$ for $y\ne0$ and $w(0)=\chi^{\prime\prime
}(0)$, is known to always converge to the unique solution of (\ref
{snrob}) regardless of the initial value $\widehat{\sigma}_{ n}^{(0)}$.
In practice, the initial value $\widehat{\sigma}_{ n}^{(0)}$ is usually
taken to be the \textsc{mad}. A discussion on the convergence of the
algorithm can be found in \citet{maro}.

For a bounded score function $\chi$, if the solution $\sigma_{\rob
}(G_0)$ of (\ref{sigmarob}) is unique, as it is the case when $\chi$ is
continuous and nondecreasing, then the functional~$\sigma_{\rob}$ is
weakly continuous at $G_0$. Weakly continuity of $\sigma_{\rob}$ at
$G_0$, that is, continuity with respect to the weak topology in
$\mathcal G$ which is given by the Prohorov metric, and consistency in
a neighborhood of $G_0$ entails robustness at $G_0$. For details, see
\citet{hub} and \citet{ham}.

%s3 #&#
\section{The estimators}\label{prop}
We consider several robust approaches in this section and define them
on a separable Hilbert space $\mathcal H$, keeping in mind that the
main application will be $\mathcal{H}=L^2(\mathcal{I})$.
From now on and throughout the paper, $\{X_i\dvtx 1 \leq i \leq n\}$ denote
realizations of the stochastic process $X\sim P$ in a separable Hilbert
space~$\mathcal{H}$. Thus, $X_{i} \sim P$ are independent stochastic
processes that follow the same law. This independence condition could
be relaxed, since we only need the strong law of large numbers to hold
in order to guarantee the results in this paper.

%s3.1 #&#
\subsection{Raw robust projection-pursuit approach}\label{rawpp}

Based on property (\ref{MAX}) of the first principal component and
given $\sigma_{\rob}(F)$ a robust scale functional, the raw (meaning
unsmoothed) robust functional principal component directions are
defined as
%
%e3.1 #&#
%
\begin{equation}\label{MAXROB}
\cases{
\displaystyle \phi_{\rob,1}(P) =\argmax_{\|\alpha\|=1}\sigma_{\rob}
(P[\alpha
]),\cr
\displaystyle \phi_{\rob,m}(P)= \argmax_{\|\alpha\|=1, \alpha\in\mathcal{B}_m}
\sigma_{\rob}(P[\alpha]),\qquad 2 \leq m,}
\end{equation}
where $P[\alpha]$ stands for the distribution of $\langle\alpha
,X\rangle
$ when $X\sim P$ and $\mathcal{B}_m=\{\alpha\in\mathcal{H}\dvtx
\langle
\alpha, \phi_{\rob,j}(P) \rangle=0, 1\le j\le m-1\}$. We will denote
the $m$th largest principal value by
%
%e3.2 #&#
%
\begin{equation}\label{lamrob}
\lambda_{\rob,m}(P)=\sigma_{\rob}^2(P[\phi_{\rob,m}]) = \max_{\|
\alpha\| =1, \alpha\in\mathcal{B}_m} \sigma^2_{\rob}(P[\alpha] ) .
\end{equation}
Since the unit ball is weakly compact, the maximum above is attained if
the scale functional $\sigma_{\rob}$ is (weakly) continuous.

Next, denote by $s^2_{n}\dvtx\mathcal{H}\to\real$ the function
$s^2_{n}(\alpha)=\sigma_{\rob}^2(P_n[\alpha])$, where\break
$\sigma
_{\rob}(P_n[\alpha])$ stands for the functional $\sigma
_{\rob}$
computed at the empirical distribu\-tion of $\langle\alpha,X_1\rangle,
\ldots, \langle\alpha,X_n\rangle$. Analogously, the mapping $\sigma
\dvtx\mathcal{H}\to\real$ stands for $\sigma(\alpha)=\sigma_{\rob
}(P[{\alpha}])$. The components in (\ref{MAXROB}) will be estimated
empirically by
%
%e3.3 #&#
%
\begin{equation}\label{estFPC}
\cases{
\wphi_{\raw,1} =\displaystyle \argmax_{\|\alpha\|=1}s_n(\alpha),\cr
\wphi_{\raw,m}= \displaystyle \argmax_{\alpha\in\widehat{\mathcal{B}}_m}
s_n(\alpha),\qquad 2 \leq m,}
\end{equation}
where $\widehat{\mathcal{B}}_m = \{ \alpha\in\mathcal{H}\dvtx \|
\alpha\|
=1, \langle\alpha, \wphi_{\raw,j} \rangle=0 , \forall 1\le
j \le
m-1\}$. The
estimators of the principal values are then computed as
%
%e3.4 #&#
%
\begin{equation}\label{esttamanno}
\wlam_{\raw,m}= s^2_{n}(\wphi_{\raw,m}),\qquad 1 \leq m .
\end{equation}

%s3.2 #&#
\subsection{Smoothed robust principal components}\label{smoothpp}

Sometimes instead of raw functional principal components, smoothed ones
are of interest. The advantages of smoothed functional PCA are well
documented; see, for instance, \citet{ris} and \citet{rasi2}. One compelling argument is that smoothing is a
regularization tool that might reveal more interpretable and
interesting features of the modes of variation for functional data. As
noted in the \hyperref[intro]{Introduction}, \citet{ris} and \citet{s} proposed two smoothing approaches by penalizing the variance and
the norm, respectively.

To be more specific, \citet{ris} estimate the first principal
component by maximizing over $\|\alpha\|=1$, the objective function
$\widehat{\var}(\langle\alpha,X\rangle)- \rho
\lceil\alpha,\alpha\rceil$, where $\widehat{\var}$ stands for the
sample variance and $\lceil \alpha
,\beta\rceil\!=\!\int_0^1\!\alpha^{\prime\prime}(t)\beta^{\prime \prime
}(t)\,dt$. Consistency for these estimators was established by
\citet{pez}.

Another regularization method proposed by \citet{s} is to penalize
the roughness through a norm defined via the penalized inner product,
$\langle\alpha,\beta
\rangle_{\tau}=\langle\alpha,\beta\rangle+\tau\lceil\alpha ,\beta
\rceil$. The smoothed first direction $\wphi_1$ is then the one that
maximizes $\widehat{\var}(\langle\alpha,X\rangle)$ over $\|
\alpha\|_{\tau}=1$. Consistency of these estimators is also established
in \citet{s} under the assumption that $\phi_j$ has finite
roughness, that is, $\lceil\phi_j,\phi_j \rceil<\infty$.

Clearly the smoothing parameters $\rho$ and $\tau$ need to converge to
$0$ in order to get consistency results.

Let us consider $\mathcal{H}_{\smooth}$, the subset of ``smooth
elements'' of $\mathcal{H}$. In order to obtain consistency results, we
will assume that $\phi_{\rob,j}(P)\in\mathcal {H}_{\smooth }$. Let
$D\dvtx\mathcal{H}_{\smooth}\,{\rightarrow}\,\mathcal{H}$~be a linear
operator, which we will refer to as the ``differentiator.'' Using $D$,
we define the symmetric positive semidefinite bilinear form \mbox{$
\lceil\cdot, \cdot\rceil\dvtx\mathcal{H}_{\smooth}\times\mathcal
{H}_{\smooth} \rightarrow\real$}, where $\lceil\alpha, \beta\rceil=
\langle D\alpha, D\beta\rangle$. The ``penalization operator'' is then
defined as $\Psi\dvtx \mathcal{H}_{\smooth} \rightarrow \mathbb{R}$,
$\Psi(\alpha) = \lceil\alpha,\alpha\rceil$, and the penalized inner
product as $\langle\alpha,\beta\rangle_{\tau }=\langle
\alpha,\beta\rangle+\tau\lceil\alpha,\beta\rceil$. Therefore,
$\|\alpha\|_{\tau}^2=\|\alpha\|^2+\tau\Psi(\alpha)$. As in
\citet{pez}, we will assume that the bilinear form is closable.
%
%re3.1 #&#
%
\begin{remark}\label{remark31}
The most common setting for functional data is to choose $\mathcal{H} =
L^2(\mathcal{I})$, $ \mathcal{H}_{\smooth} = \{ \alpha\in
L^2(\mathcal{I}), \alpha$ is twice\vspace*{1pt} differentiable, and $
\int_\mathcal{I} (\alpha^{\prime\prime }(t))^2\,dt<\infty\} $,
$D\alpha= \alpha^{\prime\prime}$ and $ \lceil \alpha, \beta\rceil=
\int_\mathcal{I} \alpha^{\prime\prime }(t)\beta ^{\prime\prime}(t) \,dt
$ so that $\Psi(\alpha) = \int_\mathcal{I} (\alpha^{\prime\prime}(t))^2
\,dt$.
\end{remark}

Let $\sigma_{\rob}(F)$ be a robust scale functional and define
$s_{n}(\alpha)$ and $\sigma(\alpha)$ as in Section~\ref{rawpp}.
Then we
can adapt the classical procedure by defining the smoothed robust
functional principal direction estimators either:

\begin{longlist}[(a)]
\item[(a)]
by penalizing the norm as
%
%e3.5 #&#
%
\begin{equation}\label{estFPCsmooth2}
\cases{
\displaystyle \wphi_{\smoothn, 1} =\argmax_{\|\alpha\|_{\tau
}=1}s_n^2(\alpha
)=\argmax_{\beta\ne0}
\frac{s_n^2(\beta)}
{\langle\beta,\beta\rangle+\tau\Psi(\beta)},%{=\argmax_{\|\beta
\vspace*{2pt}\cr
\displaystyle \wphi_{\smoothn, m}= \argmax_{\alpha\in\widehat{\mathcal{B}}_{m,
\tau,\smoothn}} s_n^2(\alpha),\qquad 2 \leq m,}
\end{equation}
where $\widehat{\mathcal{B}}_{m, \tau, \smoothn} = \{ \alpha\in
\mathcal{H}\dvtx\|\alpha\|_{\tau}=1, \langle\alpha, \wphi
_{\smoothn, j}
\rangle_{\tau}=0 , \forall 1\le j \le m-1\}$,

\item[(b)] or by penalizing the scale as
%
%e3.6 #&#
%
\begin{equation}\label{estFPCsmooth1}
\cases{\displaystyle
\wphi_{\smooths, 1}=\argmax_{\|\alpha\|=1}\{s_n^2(\alpha)-
\rho\Psi(\alpha)\},\vspace*{2pt}\cr
\displaystyle \wphi_{\smooths, m}=\argmax_{\alpha\in\widehat{\mathcal
{B}}_{m,\smooths}} \{s_n^2(\alpha)- \rho\Psi(\alpha)\},\qquad
2 \leq m,}
\end{equation}
where $\widehat{\mathcal{B}}_{m,\smooths} = \{ \alpha\in\mathcal
{H}\dvtx\|\alpha\|=1, \langle\alpha, \wphi_{\smooths, j} \rangle=0
,
\forall 1\le j \le m-1\}$.
\end{longlist}

The corresponding principal value estimators are respectively defined as
%
%e3.7 #&#
%
\begin{equation}\label{lamFPCsmooth}
\wlam_{\smooths,m} = s^2_{n}(\wphi_{\smooths,m}) \quad\mbox{and}\quad
\wlam_{\smoothn,m} = s^2_{n}(\wphi_{\smoothn,m}).
\end{equation}

%s3.3 #&#
\subsection{Sieve approach for robust functional principal components}
\label{sievepp}

Another approach, motivated by using $B$-splines as a smoothing tool,
is to consider the method of sieves. The method of sieves involves
approximating an infinite-dimensional parameter space $\Theta$ by a
sequence of finite-dimensional parameter spaces $\Theta_n$, which
depend on the sample size $n$, and then estimate the parameters on the
spaces $\Theta_n$ rather than $\Theta$.

Let $\{\delta_i\}_{i\ge1}$ be a basis of $\mathcal{H}$ and define
$\mathcal{H}_{p_n}$ as the linear space spanned by $\delta_1,\ldots,
\delta_{p_n}$ and by
$\mathcal{S}_{p_n}=\{\alpha\in\mathcal{H}_{p_n}\dvtx
\|\alpha\|=1\} $, that is,
\[
\mathcal{H}_{p_n}=\Biggl\{\alpha\in\mathcal{H}\dvtx\alpha=\sum
_{j=1}^{p_n}a_j \delta_j \Biggr\}
\]
and $
\mathcal{S}_{p_n}=\{\alpha\in\mathcal{H}\dvtx\alpha=\sum
_{j=1}^{p_n}a_j \delta_j$, such that
$\|\alpha\|^2=\sum
_{j=1}^{p_n}\sum
_{s=1}^{p_n}a_j a_s \langle\delta_j, \delta_s\rangle=1\}$.
Note that $\mathcal{S}_{p_n}$ approximates the unit sphere $\mathcal
{S}=\{\alpha\in\mathcal{H}\dvtx\|\alpha\|=1\}$. When $\{\delta_i\}
_{i\ge
1}$ is an orthonormal basis, $\|\alpha\|^2=\sum_{j=1}^{p_n}a_j^2=\ba
\trasp\ba$ where $\ba=(a_1,\ldots,a_{p_n})\trasp$, hence, the norm of
$\alpha$ equals the Euclidean norm of the vector~$\ba$. Define the
robust sieve estimators of the principal components as
%
%e3.8 #&#
%
\begin{equation}\label{sieveFPC}
\cases{\displaystyle \wphi_{\sieve, 1} =\argmax_{\alpha\in\mathcal
{S}_{p_n}}s_n(\alpha
),\vspace*{2pt}\cr
\displaystyle \wphi_{\sieve, m}= \argmax_{\alpha\in\widehat{\mathcal
{B}}_{n,m}} s_n(\alpha),\qquad 2 \leq m,}
\end{equation}
where $\widehat{\mathcal{B}}_{n,m} = \{ \alpha\in\mathcal
{S}_{p_n}\dvtx
\langle\alpha, \wphi_{\sieve,j} \rangle=0 , \forall 1\le j
\le
m-1\}$, and let the principal value estimators be
%
%e3.9 #&#
%
\begin{equation}\label{sievelamFPC}
\wlam_{\sieve,m}=s_n^2(\wphi_{\sieve,m}).\vadjust{\goodbreak}
\end{equation}
Some of the frequently used bases for functional data are the Fourier,
polynomial, spline and wavelet bases; see, for instance,
\citet{rasi2}.

To the best of our knowledge, the above sieve approach is new to
functional principal component analysis, even if one considers the
classical sieve estimators,
that is, when $\sigma_{\rob}$ in (\ref{sieveFPC}) is the standard deviation.

%s3.4 #&#
\subsection{A unified formulation for the robust projection pursuit approaches}

To help formulate a unified approach to the different estimators
considered in Sections \ref{smoothpp}, \ref{smoothpp} and \ref
{sievepp}, let the products $\rho\Psi(\alpha)$ or $\tau\Psi(\alpha
)$ be
defined as $0$ when $\rho=0$ or $\tau=0$, respectively, even when
$\alpha\notin\mathcal{H}_{\smooth}$ for which case $\Psi(\alpha
)=\infty$. Moreover, when $p_n=\infty$, define $\mathcal{H}_{p_n}=
\mathcal{H}$.
All the projection pursuit estimators considered in the previous
subsections then can be viewed as special cases of the following
general class of estimators:
%
%e3.10 #&#
%
\begin{equation}\label{estFPCgeneral}
\cases{\displaystyle
\wphi_{1}=\argmax_{\alpha\in\mathcal{H}_{p_n},\|\alpha\|
_{\tau
}=1}\{s_n^2(\alpha)- \rho\Psi(\alpha)\},
\vspace*{2pt}\cr
\displaystyle \wphi_{m}=\argmax_{\alpha\in\widehat{\mathcal{B}}_{m, \tau}}
\{s_n^2(\alpha)- \rho\Psi(\alpha)\},\qquad 2 \leq m,}
\end{equation}
where $\widehat{\mathcal{B}}_{m, \tau} = \{ \alpha\in\mathcal
{H}_{p_n}\dvtx\|\alpha\|_{\tau}=1, \langle\alpha, \wphi_{ j}
\rangle
_{\tau}=0 , \forall 1\le j \le m-1\}$.

With this\vspace*{1pt} definition and by taking $p_n=\infty$, the raw estimators
are obtained when $\rho=\tau=0$, while $ \wphi_{\smoothn, m}$ and $
\wphi_{\smooths, m}$ correspond to $\rho=0$ and $\tau=0$, respectively.
On the other hand, the sieve estimators correspond a~finite choice for
$p_n$ and $\tau=\rho=0$.

%s4 #&#
\section{Consistency results}\label{consist}

In this section, we show that under mild conditions the functionals
$\phi_{\rob,m}(P)$ and $\lambda_{\rob,m}(P)$ defined through (\ref
{MAXROB}) and (\ref{lamrob}) are weakly continuous. Moreover, we state
conditions that guarantee the consistency of the estimators defined in
Section \ref{prop}. Proofs for this section can be found in the
\hyperref[appenA]{Appendix} and in the supplemental article
[\citet{bali3a}].

To derive the consistency of the proposed estimators, we need the
following assumptions:
\begin{longlist}[(S0)]
\item[(S0)] For some $q\ge2$ and $1 \leq j \leq q$, $\phi
_{\rob,j}(P)$ are unique up to a sign change.
\item[(S1)] $\sigma\dvtx\mathcal{H}\to\real$ is a weakly
continuous function, that is, continuous with respect to the weak
topology in $\mathcal H$.
\item[(S2)] ${\sup_{\| \alpha\| = 1}}|s_n^2(\alpha) -
\sigma^2(\alpha)|\convpp0$.
\end{longlist}
Note that condition (S0) holds if and only if $\lambda_{\rob
,1}(P)> \cdots> \lambda_{\rob, q+1}(P)$.\vspace*{8pt}

\textit{Some remarks.}
(i) (S1) holds when the scale functional $\sigma_{\rob}$ is a
continuous functional (with respect to the weak topology under the
Prohorov distance).
This is because $\alpha_k$ converges weakly to $\alpha$, which implies
$\langle\alpha_k ,X\rangle\convweak\langle\alpha,X\rangle$ and
hence $\sigma_{\rob}(P[\alpha_k])\to\sigma_{\rob}(P[\alpha])$.
For the
case when the scale functional is the standard deviation,
and the underlying probability $P$ has a compact covariance operator
$\bGa_{X}$, we see from the relationship $\sigma^2(\alpha) = \langle
\alpha, \bGa_{X} \alpha\rangle$ that condition (S1) holds, even
though the standard deviation itself is not a~weakly continuous functional.

(ii) Since there exists a metric $d$ generating the weak topology in
$\mathcal H$ and that the closed ball $\mathcal{V}_r=\{\alpha\dvtx \|
\alpha
\|\le r\}$ is weakly compact, we see that (S1) implies that
$\sigma(\alpha)$ is uniformly continuous with respect to the metric $d$
and hence, with respect to the weak topology, over $\mathcal{V}_r$.

(iii) Assumption (S2) holds for the classical estimators based
on the sample variance since the empirical covariance operator,
$\widehat{\bGa}$, is consistent in the unit ball. Indeed, as shown in
\citet{dpr}, $\| \widehat{\bGa} - \bGa_{X}\| \convpp0$,
which entails that
$
{\sup_{\| \alpha\| = 1}}|s_n^2(\alpha) - \sigma^2(\alpha
)|
\leq\| \widehat{\bGa} - \bGa_{X} \| \convpp0.
$
However, this assumption may seem harder to verify for other scale
functionals since the unit sphere $\mathcal{S}=\{\|\alpha\|=1\}$ is not
compact, and $s_n^2(\alpha)$ is usually not defined through a
covariance operator estimator. To be more precise, even under some
conditions to be considered in Section \ref{fisher}, there exists a
compact operator $\bGa$ such that $\sigma(\alpha)=\langle\alpha,
\bGa
\alpha\rangle$, $s_n^2(\alpha)$ cannot be expressed as $\langle
\alpha,
\bGa_n\alpha\rangle$ for some consistent estimator $\bGa_n$ of
$\bGa$.
Corollary \ref{corollary61} in Section \ref{appen} establishes that
(S2) holds for any scale functional $\sigma_{\rob}$ that is continuous
with respect to the weak topology.

The following lemma, whose proof can be found in Section B
of the technical supplemental article given in
\citet{bali3a}, is useful for deriving the consistency and continuity of
the principal direction estimators. In this lemma and in the subsequent
theorems, it should be noted that $\langle\wphi,\phi\rangle^2
\rightarrow1$ implies,\vspace*{1pt} under the same mode of convergence, that the
sign of $\wphi$ can be chosen so that $\wphi\rightarrow\phi$.

For the sake of simplicity, denote by $\lambda_{\rob,j}=\lambda
_{\rob
,j}(P)$ and $\phi_{\rob,j}=\phi_{\rob,j}(P)$.
%
%le4.1 #&#
%
\begin{lemma}\label{lemma41}
Let $\wphi_m\in\mathcal
{S}$ be such that $\langle\wphi_m, \wphi_j\rangle\stackrel{\mathit{a.s.}}{\longrightarrow}0$ for
$j\ne
m$ and assume that \textup{(S0)} and \textup{(S1)} hold. Then:
\begin{longlist}[(b)]
\item[(a)] If $\sigma^2(\wphi_1)\stackrel{\mathit{a.s.}}{\longrightarrow}\sigma^2(\phi_{\rob,1})$, then
$\langle\wphi_1,\phi_{\rob,1}\rangle^2\stackrel{\mathit{a.s.}}{\longrightarrow}1$.
\item[(b)] Given $2 \leq m \leq q$, if $\sigma^2(\wphi_m)\stackrel{\mathit{a.s.}}{\longrightarrow}
\sigma
^2(\phi_{\rob,m})$ and $\wphi_s \stackrel{\mathit{a.s.}}{\longrightarrow}\phi_{\rob,s}$, for $1\le
s\le
m-1$, then $\langle\wphi_m,\phi_{\rob,m}\rangle^2\stackrel{\mathit{a.s.}}{\longrightarrow}1$.
\end{longlist}
\end{lemma}

Let $d_{\pr}(P,Q)$ stand for the Prohorov distance between the
probability measures $P$ and $Q$. Thus, $P_n \convweak P$ if and only
if $d_{\pr}(P_n, P) \to0$. Theorem~\ref{theorem41} below establishes
the continuity of the functionals defined in (\ref{MAXROB}) and~(\ref
{lamrob}), and hence the asymptotic robustness of the estimators
derived from them, as defined in \citet{ham}. This can be seen just
by replacing almost sure convergence by convergence in its statement.
As it will be shown in Section \ref{appen}, the uniform convergence
required in assumption (ii) below is satisfied, for instance, if
$\sigma
_{\rob}$ is a continuous scale functional when $P_n \convweak P$.

To accommodate data driven smoothing parameters a more general
framework is considered in Theorem \ref{theorem41}, which allows for
the smoothing parameters $\tau_n$ and $\rho_n$ to be random, such that
$\tau_n\convpp0$ and $\rho_n\convpp0$.
%
%th4.1 #&#
%
\begin{theorem}\label{theorem41}
Let $P_n$ be a sequence
of probability measures and $\tau=\tau_n\ge0$, $\rho=\rho_n\ge0$ be
random smoothing parameters. Denote\vspace*{1pt} by $\sigma_n^2(\alpha)=\sigma
_{\rob
}^2(P_{n}[\alpha])$ and define $ \wlam_{ m}=\sigma_n^2(\wphi_{m})$ with
%
%e4.1 #&#
%
\begin{equation}\label{estFPC2}
\cases{\displaystyle \wphi_{1} =\argmax_{\|\alpha\|_{\tau}=1} \{
\sigma_n^2(\alpha
)- \rho\Psi(\alpha)\},
\vspace*{2pt}\cr
\displaystyle \wphi_{ m}=\argmax_{\alpha\in\widehat{\mathcal{B}}_{m, \tau}}
\{\sigma_n^2(\alpha)- \rho\Psi(\alpha)\},\qquad 2 \leq m,}
\end{equation}
where $\widehat{\mathcal{B}}_{m, \tau} = \{ \alpha\in\mathcal
{H}\dvtx \|
\alpha\|_{\tau}=1, \langle\alpha, \wphi_{j} \rangle_{\tau}=0
,
\forall 1\le j \le m-1\}$. Let $P$ be a~probability measure
satisfying \textup{(S0)}. Assume that:
\begin{longlist}
\item \textup{(S1)} holds;
\item ${\sup_{\|\alpha\|=1}}|\sigma_n^2(\alpha) - \sigma
_{\rob
}^2(P[\alpha])|\convppp0$;\vspace*{1pt}
\item $\tau_n\convppp0$ and $ \rho_n \convppp0$;
\item moreover, if $\tau_n>0 $ or $ \rho_n> 0 $, for all $n\ge
n_0$, assume that $\phi_{\rob,j} \in\mathcal{H}_{\smooth}$, that is,
$\Psi( \phi_{\rob,j})<\infty$, for all $ 1 \le j\le q$.
\end{longlist}
Then:
\begin{longlist}[(d)]
\item[(a)] $ \wlam_{1}\convppp\lambda_{\rob,1}$ and $\sigma
^2(\wphi
_{1})\convppp\sigma^2(\phi_{\rob,1})$. Moreover, $\rho\Psi(\wphi_{1})
\convppp0$ and $\tau\lceil\wphi_{1},\wphi_{1} \rceil\convppp0$, and
so $\|\wphi_{1}\| \convppp1$;\vspace*{1pt}
\item[(b)] $\langle\wphi_{1},\phi_{\rob,1}\rangle^2\convppp1$;
\item[(c)] for any $2 \leq m \leq q$, if $\wphi_{\ell} \convppp\phi
_{\rob
,\ell}$, $\tau\Psi(\wphi_{ \ell}) \convppp0$ and $\rho\Psi(\wphi_{
\ell}) \convppp0$ for $1\le\ell\le m-1$, then $ \wlam_{m}\convppp{
\sigma^2(\phi_{\rob,m})}$ and $\sigma^2(\wphi_{m})\convppp\sigma
^2(\phi
_{\rob,m})$. Moreover, $\rho\Psi(\wphi_{m}) \convppp0$, $\tau\Psi
(\wphi_{ m}) \convppp0$ and so, $\|\wphi_{m}\| \convppp1$;
\item[(d)] for $1\leq m \leq q$, $\langle\wphi_{m},\phi_{\rob
,m}\rangle
^2\convppp1$.
\end{longlist}
\end{theorem}

Note that assumption (ii) corresponds to (S2) when $P_n$ is
the empirical probability measure. On the other hand, when $\sigma
_{\rob
}(\cdot)$ is a continuous scale functional, Theorem \ref{theorem62}
implies that (ii) holds whenever \mbox{$d_{\pr}(P_n,P)\convpp0$}. Moreover,
if $\sigma_{\rob}(\cdot)$ is a continuous scale functional and $P$
satisfies (S0), Theorem \ref{theorem41} entails the continuity
of the functionals $\phi_{\rob,j}(\cdot)$ and $\lambda_{\rob
,j}(\cdot)$
at~$P$, for $1\le j\le q$, and so the proposed estimators are
qualitatively robust and consistent. In particular, the estimators are
robust at any elliptical distribution $\mathcal{E}(\mu, \bGa)$, as
defined in Section \ref{fisher}, such that the largest $q+1$
eigenvalues of the operator $\bGa$ are all distinct.

Theorem \ref{theorem41} establishes the consistency of the raw
estimators of the principal components under (S0) to (S2)
by taking $\rho=\tau=0$. It also shows that proposals (\ref
{estFPCsmooth2}) and (\ref{estFPCsmooth1}) give consistent estimators
if $\phi_{\rob,j} \in\mathcal{H}_{\smooth}$, $1 \le j \le q$.
In\vspace*{1pt}
\citet{bali2}, it is shown that the estimators $\wphi_{\smoothn,m}$ and
$\wlam_{\smoothn,m}$ defined in (\ref{estFPCsmooth2}) and (\ref
{lamFPCsmooth}) are still consistent if $\phi_{\rob,j} \in\overline
{\mathcal{H}}_{\smooth}$, $1 \le j \le q$, where $\overline{\mathcal
{H}}_{\smooth}$ stands for the closure of $\mathcal{H}_{\smooth}$. The
condition $\phi_{\rob,j} \in\overline{\mathcal{H}}_{\smooth}$
generalizes the assumption of smoothness, $\phi_{\rob,j} \in\mathcal
{H}_{\smooth}$, required in \citet{s} and holds, for example,
when $\mathcal{H}_{\smooth}$ is a dense subset of $\mathcal{H}$.

Theorem \ref{theorem42} establishes the consistency of the estimators
of the principal directions defined through the sieve approach given in
(\ref{sieveFPC}). Below we give a
separate statement for the consistency of the sieve estimators to avoid
imposing additional burdensome assumptions, such as smoothness
conditions for the basis elements, whenever either
$\tau\ne0$ or $\rho\ne0$ in (\ref{estFPCgeneral}). Its proof is
relegated to the Section C of the technical supplement
[\citet{bali3a}].
%
%th4.2 #&#
%
\begin{theorem}\label{theorem42}
Let $\wphi_{\sieve,m}$
and $\wlam_{\sieve,m}$ be the estimators defined in (\ref{sieveFPC})
and (\ref{sievelamFPC}), respectively. Under \textup{(S0)} to
\textup{(S2)}, if $p_n\to\infty$, then:
\begin{longlist}[(c)]
\item[(a)] $ \wlam_{\sieve,1}\convppp{ \sigma^2(\phi_{\rob,1})}$ and
$\sigma^2(\wphi_{\sieve,1})\convppp\sigma^2(\phi_{\rob,1})$.
\item[(b)] Given $2 \leq m \leq q$, if $\wphi_{\sieve,\ell} \convppp
\phi
_{\rob,\ell}$, for $1\le\ell\le m-1$, then $ \wlam_{\sieve
,m}\convppp{
\sigma^2(\phi_{\rob, m})}$ and $\sigma^2(\wphi_{\sieve,m})\convppp
\sigma
^2(\phi_{\rob,m})$.
\item[(c)] For $1\leq m \leq q$, $\langle\wphi_{\sieve,m},\phi
_{\rob
,m}\rangle^2\convppp1$.
\end{longlist}
\end{theorem}

%s5 #&#
\section{Fisher-consistency under elliptical distributions}\label{fisher}
The results in Section~\ref{consist} ensure that, under mild
conditions, the estimates of the principal directions converge almost
surely to $\phi_{\rob,m}$ defined in (\ref{MAXROB}). An important point
to highlight is what the functions $\phi_{\rob,m}$ represent, at least
in some particular situations. This section focuses on showing that,
for the functional elliptical families defined in \citet{bali},
the functionals~$\phi_{\rob,m}(P)$ and~$\lambda_{\rob,m}(P)$
have a simple interpretation. In particular, our results hold for the
functional elliptical family, but are not restricted to it. We recall
here their definition for the sake of completeness.

Let $X$ be a random element in a separable Hilbert space $\mathcal H$
and $\mu\in\mathcal{H}$. Let $\bGa\dvtx\mathcal{H} \rightarrow
\mathcal
{H}$ be a self-adjoint, positive semidefinite and compact operator. We
will say that $X$ has an elliptical distribution with parameters $(\mu,
\bGa)$, denoted as $X\sim\mathcal{E}(\mu, \bGa)$, if for any linear
and bounded operator $A\dvtx\mathcal{H} \rightarrow\real^d$, $AX$~has a
multivariate elliptical distribution with parameters $A\mu$ and $A\bGa
A^* $, that is, $AX\sim\mathcal{E}_{d}(A\mu, A\bGa A^*)$, where
$A^*\dvtx
\mathbb{R}^p\rightarrow\mathcal{H}$ stands for the adjoint operator of~%
$A$. As in the finite-dimensional setting, if the covariance operator,
$\bGa_X$, of~$X$ exists, then $\bGa_X = {a} \bGa$, for some $a\in
\real$.

The elliptical distributions in $\mathcal{H}$ include the Gaussian
distributions. Other elliptical distributions can
be obtained from the following construction. Let~$V_1$ be a Gaussian
element in $\mathcal{H}$ with zero mean and covariance operator $\bGa
_{V_1}$, and
let $Z$ be a random variable independent of $V_1$. Given $\mu\in
\mathcal{H}$, define $X=\mu+Z V_1$. Then, $X$ has an elliptical
distribution $\mathcal{E}(\mu, \bGa)$
with the operator $\bGa$ being proportional to $\bGa_{V_1}$. Note that
$\bGa$ exist even if the second moment of~$X$ do not exist.
For random elements which admit a finite Karhunen--Lo\`eve expansion,
that is, $X(t) = \mu(t) + \sum_{j=1}^q \lambda_j^{1/2} U_j \phi_j(t)$,
the assumption that $X$ has an elliptical distribution is analogous to
assuming that $\bU=(U_1,\ldots,U_q)\trasp$ has a spherical distribution.

Lemma \ref{lemma51} below states the Fisher-consistency of the
functionals defined through (\ref{MAXROB}) under the following
assumption:
\begin{longlist}[(S3)]
\item[(S3)] There exists a constant $c>0$ and a self-adjoint,
positive semidefinite and compact operator $\bGa_0$, such that for any
$\alpha\in\mathcal{H}$, we have $\sigma^2(\alpha)={c}\langle
\alpha,
\bGa_0\alpha\rangle$.
\end{longlist}
Its proof follows immediately and is thus omitted. Note that
(S3) entails that the function $\sigma\dvtx \mathcal{H} \rightarrow\real$
defined as $\sigma(\alpha) = \sigma_{\rob}(P[\alpha] )$ is weakly
continuous, hence (S1) holds. Besides, as a consequence of Lemma
\ref{lemma51}, (S0) holds under~(S3) if the $q$ largest
eigenvalues of $\bGa_0$ are distinct.
%
%le5.1 #&#
%
\begin{lemma}\label{lemma51}
Let $\phi_{\rob,m}$ and $\lambda_{\rob,m}$ be the functionals defined
in (\ref{MAXROB}) and~(\ref{lamrob}), respectively. Let $X\sim P$ be a
random element such that \textup{(S3)} holds. Denote by
$\lambda_1\ge\lambda_2\ge\cdots$ the eigenvalues of $\bGa_0$ and by
$\phi_j$ the eigenfunction of $\bGa_0$ associated to~$\lambda_j$.
Assume that for some $q\ge2$, and for all $1\le j\le q$, $\lambda
_1>\lambda_2>\cdots>\lambda_q>\lambda_{q+1}$. Then, we have that
$\phi
_{\rob,j}(P) = \phi_j$ and $\lambda_{\rob,j}(P) = c \lambda_j$.
\end{lemma}

For any distribution possessing finite second moments, if the scale
functional is taken to be the standard deviation, then (S3)
holds with $\bGa_0=\bGa_X$. When considering a robust scale functional,
(S3) holds if $X$ has an elliptical distribution $\mathcal{E}(\mu
, \bGa)$ taking $\bGa_0=\bGa$, and so Lemma \ref{lemma51} entails that
the functionals $\phi_{\rob,j}(P)$ defined through (\ref{MAXROB}) are
Fisher-consistent. As in the finite-dimensional setting, the scale
functional $\sigma_{\rob}$ can be calibrated to attain
Fisher-consistency of the principal values.

Assumption (S3) ensures that we are estimating the target
directions. It may seem restrictive since it is difficult to verify outside
the family of elliptical distributions except when the scale is taken
to be the standard deviation. However, even in the finite-dimensional case,
asymptotic properties have been derived only under similar restrictions
when considering projection-pursuit estimators. For instance, both
\citet{li} and \citet{crouruiz2} assume an underlying
elliptical distribution in order to obtain consistency results and
influence functions, respectively. Also, in \citet{cui} the
influence function of the projected data is assumed to be of the form
$h(\bx,\ba)=2 \sigma(F[\ba])\operatorname{IF}(\bx, \sigma_{\ba} ;F_0)$, where
$F[\ba]$ stands for the distribution of $\ba\trasp\bx$ when $\bx
\sim
F$. This condition, though, primarily holds when the distribution is elliptical.
%
%re5.1 #&#
%
\begin{remark}\label{remark51}
An alternative to the robust
projection pursuit approach for functional principal components is to
consider the spectral value decomposition of a robust covariance or
scatter operator. The spherical principal components, noted in the
\hyperref[intro]{Introduction}, which were proposed by \citet{loc} and
further developed by \citet{ger08}, apply this approach using the
spatial covariance operator. The spatial covariance operator is defined
to be
\[
\bV= \esp\bigl({(X - \eta) \otimes(X- \eta)}/{\|X - \eta\|^2}\bigr)
\]
with $\eta$ being the spatial median, that is,
%
%e5.1 #&#
%
\begin{equation}\label{spatialm}
\eta= \argmin_{\theta\in\mathcal{H}}\esp( \|X - \theta\| - \| X \| ).
\end{equation}
The spatial median is sometimes referred to as the multivariate $L^1$
median, but this is a misnomer since the the norm in (\ref{spatialm})
is the $L^2$ norm. Note that when the norm is replaced by the square of
the norm in (\ref{spatialm}), the resulting parameter is the mean.

\citet{ger08} proved the Fisher-consistency of the eigenfunctions of
the spatial covariance operator, but under the additional assumption
that $X$ has a finite Karhunen--Lo\`eve expansion so that $\bV$ has
only a finite number of nonzero eigenvalues, which is
essentially the multivariate setting. Unlike the projection pursuit
approach, though, under an elliptical model the eigenvalues of $\bV$
are not proportional to the eigenvalues of the shape parameter~$\bGa$.
Consequently, as discussed, for example, in \citet{mar}, \citet{bofr}
and \citet{vis}, this implies that even if the
second moments exist, the amount of variance explained by a~principal
component variable is not equivalent to the ratio of the eigenvalue to
the sum of all the eigenvalues. That is, $\wtlam_k/\sum_{j=1}^{\infty}
\wtlam_j$ is not the same as the explained proportion $\lambda_k/\sum
_{j=1}^{\infty} \lambda_j$, where
$\wtlam_k$ and $\lambda_k$ are the $k$th largest eigenvalue of~$\bV$
and $\bGa$ respectively.

In the multivariate setting, it is also known that the eigenvectors
obtained from the sample spatial covariance matrix can be extremely
inefficient estimates whenever the eigenvalues differ greatly; see
\citet{croux99}. Intuitively, the reason for this inefficiency is that the
spatial covariance matrix down-weights observations according to their
Euclidean distance from the center. This seems reasonable when the
distribution is close to being spherical, but less so when there are
strong dependencies in the variables. In some sense, this is the
antithesis of PCA, since one is usually interested in principal
components when one suspects the latter.

As noted in \citet{maro}, there is a vast literature on
robust estimates for covariance matrices, such as $M$-estimates,
$S$-estimates and the $\mathit{MCD}$, among others. These estimates downweight
observations relative to the shape of the data cloud. It may seem
reasonable then to try to extend these estimates to the functional
setting. An important feature of these estimates, though, is that they
are affine equivariant, and as shown in \citet{tyler}, this implies
that, when the sample size is no greater than the dimension plus one,
such estimates are simply proportional to the sample covariance matrix.
In the functional data setting, the sample size is always less than the
dimension, which is infinite. Thus, at this time, we view the robust
projection-pursuit approach as more viable.
\end{remark}

%s6 #&#
\section{Some uniform convergence results}\label{appen}
In this section, we show that when the scale functional is continuous
with respect to the Prohorov distance, (S2) and more generally,
condition (ii) in Theorem \ref{theorem41} hold whenever $P_n\convweak
P$. To derive these results, we will first state some properties
regarding the weak convergence of empirical measures that hold not only
in $L^2(\mathcal{I})$ but in any complete and separable metric space.
These properties may be useful in other settings. The proofs for the
theorems in this section are relegated to Section D of the
supplemental article [\citet{bali3a}].

Let $\mathcal M$ be a complete and separable metric space (Polish
space) and $\mathcal{B}$ the Borel $\sigma$-algebra of $\mathcal M$.
Lemma \ref{lemma61}, which is a restatement of Theorem 3 in
\citet{varad}, ensures that the empirical measures converge
weakly almost surely on a Polish space to the probability measure
generating the observations.
%
%le6.1 #&#
%
\begin{lemma}\label{lemma61}
Let $(\Omega,\mathcal{A},\prob)$ be a probability space and
$X_n\dvtx\Omega
\rightarrow\mathcal{M}$, $n \in\mathbb{N}$, be a sequence of
independent and identically distributed random elements such that
$X_i\sim P$. Assume that $\mathcal M$ is a Polish space, and denote by
$P_n$ the the empirical probability measure, that is, $P_n(A)=
({1}/{n}) \sum_{i=1}^n I_A(X_i)$ with $I_A(X_i)=1$ if $X_i\in A$ and 0
elsewhere. Then, $P_n \convweak P$ almost surely, that is, $d_{\pr
}(P_n,P)\convppp0$.
\end{lemma}

Let $P$ be a probability measure in $\mathcal M$, a separable Banach
space, and let $\mathcal{M}^*$ denote the dual space. For a given $f
\in\mathcal{M}^*$, define $P[f]$ as the real measure of the random
variable $f(X)$, with $X \sim P$.
%
%th6.1 #&#
%
\begin{theorem}\label{theorem61}
Let $\{P_n\}_{n \in\mathbb{N}}$ and $P$ be probability measures
defined on~$\mathcal M$ such that $d_{\pr}(P_n, P) \to0$. Then, $\sup
_{\| f \|_{*} = 1} d_{\pr}(P_n[f], P[f]) \to0$.
\end{theorem}

When the Banach space $\mathcal M$ above is a separable Hilbert space
$\mathcal{H}$, the Riesz representation theorem implies that for $f
\in
\mathcal{H}^*$ with $\|f \|_* =1$, there exists $\alpha\in\mathcal
{H}$ such that $f(X)=\langle\alpha, X\rangle$. The following result
states that when~$\sigma_{\rob}$ is a continuous scale functional,
uniform convergence can be attained and so, assumption (ii) in Theorem
\ref{theorem41} is satisfied.
%
%th6.2 #&#
%
\begin{theorem}\label{theorem62}
Let $\{P_n\}_{n \in
\mathbb{N}}$ and $P$ be probability measures defined on a~separable
Hilbert space $\mathcal{H }$, such that $d_{\pr}(P_n, P) \to0$. Let
$\sigma_{\rob}$ be a continuous scale functional. Then,
$\sup_{\| \alpha\| = 1 } |\sigma_{\rob}(P_n[\alpha]) - \sigma
_{\rob
}(P[\alpha])| \longrightarrow0$.
\end{theorem}

Using Lemma \ref{lemma61} and Theorem \ref{theorem62}, we get the
following result that shows that (S2) holds if $\sigma_{\rob}$
is a continuous scale functional.
%
%co6.1 #&#
%
\begin{corollary}\label{corollary61}
Let $P$ be a probability measure in a separable Hilbert space $\mathcal
{H }$, $P_n$ be the empirical measure of a random sample $X_1, \ldots,
X_n$ with $X_i \sim P$, and $\sigma_{\rob}$ be a continuous scale
functional. Then we have that\break
$
{\sup_{\| \alpha\| = 1 }} |\sigma_{\rob}(P_n[\alpha]) - \sigma
_{\rob
}(P[\alpha])| \convppp0.
$
\end{corollary}

%s7 #&#
\section{Selection of the smoothing parameters}\label{smoothpar}
The selection of the smoothing parameters is an important practical
issue. The most popular general approach to address such a selection problem
is to use the cross-validation methods. In nonparametric regression,
the sensitivity of $L^2$ cross-validation methods to outliers
has been pointed out by \citet{wang} and by \citet{can}, among others. The latter also proposed more robust
alternatives to
$L^2$ cross-validation. The idea of robust cross-validation can be
adapted to the present situation. Assume for the moment that we are
interested in
a fixed number, $\ell$, of components. We propose to proceed as follows:

\begin{longlist}[(CV1)]
\item[(CV1)]
Center the data. That is, define $\widetilde
{X}_{i}=X_{i}-\widehat
{\mu}$ where $\widehat{\mu}$ is a robust location estimator, such as
the functional spatial median defined in \citet{ger08}.
\item[(CV2)] For the penalized roughness approaches and for each $m$ in the
range $1\le m\le\ell$ and ${\tau> 0}$, let $\wphi_{m,\tau}^{(-j)}$
denote the robust estimator of the $m$th principal direction computed
without the $j$th observation.
\item[(CV3)] Define $ X_{j}^{\bot}(\tau)= \widetilde{X}_{j}-\pi(\widetilde
{X}_{j};{\widehat{\mathcal{L}}_\ell^{(-j)}})$, for $1\le j \le n$,
where $\pi(X;\mathcal{L})$ stands for the orthogonal projection of $X$
onto the linear (closed) space $\mathcal{L}$, and $\widehat{\mathcal
{L}}_\ell^{(-j)}$ stands for the linear space spanned by $\wphi
_{1,\tau
}^{(-j)}, \ldots, \wphi_{\ell,\tau}^{(-j)}$.
\item[(CV4)] Let $\mathrm{RCV}_\ell(\tau)= \sigma_n^2( \|X_{1}^{\bot}(\tau)\|,
\ldots, \|
X_{n}^{\bot}(\tau)\| )$, where $\sigma_n$ is a robust measure of scale
about zero. We then choose $\tau_n$ to be the value of $\tau$ which
minimizes $\mathrm{RCV}_\ell(\tau)$.
\end{longlist}

By a robust measure of scale about zero, we mean that no location
estimator is applied to center the data. For instance, in the classical
setting, one takes $\sigma_n^2(z_1,\ldots,z_n)=(1/n)\sum_{i=1}^n z_i^2$,
while in the robust situation, one might choose $\sigma_n(z_1,\ldots
,z_n)=\median(|z_1|,\ldots,|z_n|)$ or to be
an $M$-estimator satisfying equation (\ref{snrob}) when setting
$\widehat{\mu}_n =0$.

For large sample sizes, it is well understood that cross-validation
methods can be computationally prohibitive. In such cases, $K$-fold
cross-validation
provides a useful alternative. In the following, we briefly describe a
robust $K$-fold cross-validation procedure suitable for our proposed estimates.

\begin{longlist}[(K1)]
\item[(K1)]
First center the data as above, using $\widetilde
{X}_{i}=X_{i}-\widehat{\mu}$.

\item[(K2)] Partition the centered data set $\{\widetilde{X}_{i}\}$ randomly
into $K$ disjoint subsets of approximately equal sizes with the $j$th
subset having size \mbox{$n_j\ge2$}, $\sum_{j=1}^K n_j =n$. Let $\{
\widetilde
{X}_i^{(j)}\}_{1\le i\le n_j}$ be the elements of the $j$th subset,
and\break
$\{\widetilde{X}_i^{(-j)}\}_{1\le i\le n-n_j}$ denote the elements in
the complement of the $j$th subset. The set $\{\widetilde{X}_i^{(-j)}\}
_{1\le i\le n-n_j}$ will be the training set and $\{\widetilde
{X}_i^{(j)}\}_{1\le i\le n_j}$ the validation set.\vspace*{1pt}

\item[(K3)] Similar to step (CV2) but with leaving the $j$th validation subset
$\{
\widetilde{X}_i^{(j)}\}_{1\le i\le n_j}$ out instead of the $j$th observation.

\item[(K4)] Define $ X_{j}^{(j)\bot}(\tau)$ the same way as in step (CV3), using
the validation set. For instance, $ X_{i}^{(j)\bot}(\tau)= \widetilde
{X}_i^{(j)}-\pi(\widetilde{X}_i^{(j)};{\widehat{\mathcal{L}}_\ell
^{(-j)}})$, $1\le i\le n_j$, where $\widehat{\mathcal{L}}_\ell^{(-j)}$
stands for the linear space spanned by $\wphi_{1,\tau}^{(-j)}, \ldots,
\wphi_{\ell,\tau}^{(-j)}$.

\item[(K5)] Let\vspace*{1pt} $\mathrm{RCV}_{\ell, \KCV}(\tau)= \sum_{j=1}^K
\sigma_n^2( \| X_{1}^{(j)\bot}(\tau)\|,\ldots,
\|X_{n_j}^{(j)\bot}(\tau)\| )$, and choo\-se~$\tau_n$ to be the value of
$\tau$ which minimizes $\mathrm{RCV}_{\ell, \KCV}(\tau)$.
\end{longlist}

A similar approach can be developed to choose $p_n$ for the sieve
estimators.

%s8 #&#
\section{Monte Carlo study}\label{monte}

The results of Section \ref{consist} established under general
conditions the consistency of the various robust projection pursuit
approaches to functional principal components analysis. The classical
approach to functional principal components analysis also yields
consistent estimates, provided second moment exists. A study of the
influence functions and the asymptotic distributions of the various
procedures would be useful to compare them. We leave these important
and challenging theoretical problems for future research. For now, to
help illuminate possible differences in the various approaches, we
present in this section the results of a Monte Carlo study.

%s8.1 #&#
\subsection{Algorithms}\label{algo}

All the computational methods to be considered here are modifications
of the basic \textsc{cr} algorithm proposed by \citet{crouruiz} for the computation of principal components using
projection-pursuit. The basic algorithm applies to multivariate data,
say $m$-dimensional, and requires a search over projections in $\real
^m$. The \textsc{grid} algorithm described in \citet{crfil} can\vadjust{\goodbreak}
also be considered, in particular, when the number of variables $m$ is
larger than the sample size $n$. For the sake of completeness, we
briefly recall the \textsc{cr} algorithm.

Let $\bY= ({\mathbf{y}}_1, \ldots, {\mathbf{y}}_n)$ be the sample in $\real
^m$, and let
$\widehat{\bmu}_n(\bY)$ be a location estimate computed from this
sample. Let $1\le q\le m$ be the desired number of components to be
computed and
denote by $\xi_n$ the univariate projection index to be maximized. In
the multivariate setting, the index $\xi_n$ corresponds to a robust
univariate scale statistic.

\begin{longlist}[(CR1)]
\item[(CR1)]
For $k = 1$, set ${\mathbf{y}}^{(1)}_i = {\mathbf{y}}_i- \widehat{\bmu
}_n(\bY)$. Let
the set of candidate directions for the first
principal direction be $\mathcal{A}_{n,1}(\bY) = \{{\mathbf{y}}^{(1)}_i/\nu_i,
1\le i\le n\}$ where $\nu_i^2={{\mathbf{y}}^{(1)}_i}\trasp{\mathbf{y}}^{(1)}_i$. We then define
${\mathbf{v}}_{1} = \arg\max_{\ba\in\mathcal{A}_{n,1}(\bY)} \xi
_n(\ba\trasp
{\mathbf{y}}_1,\ldots, \ba\trasp{\mathbf{y}}_n)$.

\item[(CR2)] For $2\le k\le q$, define recursively $z^{(k-1)}_i={
\mathbf{v}}_{k-1}\trasp
{\mathbf{y}}_i$ and ${\mathbf{y}}^{(k)}_i =\break{\mathbf{y}}^{(k-1)}_i- z^{(k-1)}_i
{\mathbf{v}}_{k-1} ={\mathbf{y}}
^{(1)}_i-\pi_{\mathcal{V}_{k-1}}({\mathbf{y}}^{(1)}_i)$,
where $\pi_{\mathcal{V}_{k-1}}({\mathbf{y}})$ stands for the orthogonal
projection of ${\mathbf{y}}$ over the linear space $\mathcal{V}_{k-1}$ spanned
by ${\mathbf{v}}_{1},\ldots, {\mathbf{v}}_{k-1}$. Let the set of candidate
directions for
the $k$th principal direction be\break
$\mathcal{A}_{n,k}(\bY) = \{{\mathbf{y}}^{(k)}_i/\nu_i, 1\le i\le n\}
$ where
$\nu_i^2={{\mathbf{y}}^{(1k)}_i}\trasp{\mathbf{y}}^{(k)}_i$, and define
${\mathbf{v}}_{k} =\break\arg\max_{\ba\in\mathcal{A}_{n,k}(\bY)} \xi
_n(\ba\trasp
{\mathbf{y}}_1,\ldots, \ba\trasp{\mathbf{y}}_n)$.
\end{longlist}
The vectors ${\mathbf{v}}_{k}$ then yield approximations to the $k$th
principal direction, for $1\le k \le q$, and approximate scores for the
$k$th principal variable are given by $z^{(k)}_i={\mathbf{v}}_{k}\trasp
{\mathbf{y}}
_i$, for $1\le i\le n$. As mentioned in \citet{crouruiz},
the \textsc{cr} algorithm makes no smoothness assumptions on the index
$\xi_n$, is simple and fast, and requires only $O(n)$ computing space.

To apply the algorithm to functional data when considering either the
raw or a penalized approach, we first discretize the domain of the
observed function over $m$ equally spaced points in $\mathcal{I}=[
-1,1]$. The resulting multivariate observations are then ${\mathbf{y}}
_i=(X_i(t_1),\ldots, X_i(t_m))\trasp$, where $t_0=-1<t_1<\cdots
<t_m<t_{m+1}=1$. The index $\xi_n$ in the algorithm depends on the
approach being used. For instance, for the raw robust estimate and for
those penalizing the norm the index is a robust scale. On the other
hand, for the robust penalized scale approach the index is the square
of the robust scale plus the penalization term. Also, for the penalized
norm approach the orthogonal projection $\pi_{\mathcal{V}_{k-1}}({\mathbf{y}})$
in step (CR2) is with respect to the inner product $\langle\cdot, \cdot
\rangle_{\tau}$ so, the finite-dimensional inner product is modified as
in \citet{s}. The resulting directions ${\mathbf{v}}_k$ then give
numerical approximations for $\{\widehat{\phi}_k(t_1),\ldots,
\widehat
{\phi}_k(t_m) \}$. One can then interpolate or use smoothing methods to
obtain $\widehat{\phi}_k(t)$ for
$t \in\mathcal{I}$.

For the sieve approach, let $\widetilde{\delta}_1,\ldots, \widetilde
{\delta}_{p_n}$ be an orthonormal basis for $\mathcal{H}_{p_n}$, which
can be obtained by using Gramm--Schmidt on the original basis. For
$\alpha\in\mathcal{H}_{p_n}$, we then have $\langle\alpha, X_i
\rangle= \ba\trasp{\mathbf{y}}_i $, where $\alpha= \sum_{j=1}^{p_n} a_j
\widetilde{\delta}_j$, $\ba= (a_1,\ldots,\allowbreak a_{p_n})\trasp$, ${\mathbf{y}}
_i=(\langle X_i, \widetilde{\delta}_1\rangle,\ldots, \langle X_i,
\widetilde{\delta}_{p_n}\rangle)\trasp$. Consequently, we can take
$m=p_n$ and apply the \textsc{cr} algorithm to the inner scores ${\mathbf{y}}
_i$. The inner scores are computed numerically by approximating the
integrals over a grid of $50$ points.\vspace*{1pt} A numerical approximation for
$\widehat{\phi}_k(t)$ is then given by $\sum_{j=1}^{p_n} v_{k,j}
\widetilde{\delta}_j(t)$ with ${\mathbf{v}}_k=(v_{k,1},\ldots
,v_{k,p_n})\trasp$.

%s8.2 #&#
\subsection{The estimators}\label{estimators}
There are three main characteristics that distinguish the different
estimators: the method of
centering in the first step of the \textsc{cr} algorithm, the scale
function being used and the type of smoothing method.\vspace*{8pt}

\textit{Centering}: For classical procedures, that is, those based on
the standard deviation, we use the sample mean as the centering point
for the trajectories. For the robust procedures, that is, those based
on \textsc{mad} or \textsc{$M$-scale}, we center the data by using
either (i) the component-wise sample median, that is, the median at
each time point, or (ii) the sample spatial median; see~(\ref{spatialm}).
It turns out that the two robust centering methods produced similar
results, so only the results for the spatial median are reported.\vspace*{8pt}

\textit{Scale function}: Three scale estimators are considered here:
the classical standard deviation (\textsc{sd}), the median absolute
deviation (\textsc{mad}) and an $M$-estimator of scale (\textsc
{$M$-scale}). For the $M$-estimator, we use the score function (\ref
{funcionBT}) introduced by \citet{beatuk}, as discussed in
Section~\ref{scale}, with $c=1.56$, $\delta= 1/2$.\vspace*{8pt}

\textit{Smoothing parameters $\tau$ and $\rho$}: For both the classical
and robust procedures defined in Section \ref{smoothpp}, a penalization
depending on the $L^2$ norm of the second derivative, multiplied by a
smoothing factor, is included, that is, $\Psi(\alpha)=\int_{-1}^1
(\alpha^{\prime\prime}(t))^2\,dt$. Again the integral is computed over
the same grid of points $t_1,\ldots, t_m$, and the second derivative of
$\alpha$ at $t_i$ is approximated by $\{\alpha(t_{i+1})-2
\alpha
(t_{i})+\alpha(t_{i-1}) \}/(t_{i+1}-t_i)^2$, since we choose an
equidistant grid of points. Note that when $\rho=\tau=0$, the raw
estimators defined in Section~\ref{rawpp} are obtained.\vspace*{8pt}

\textit{Sieve}: Two different sieve basis are considered: the Fourier
basis obtained by taking $\delta_j$ to be the Fourier basis functions,
and the cubic $B$-spline basis functions. The Fourier basis used in the
sieve method is the same basis used to generate the data.\vspace*{8pt}

In all tables, the estimators corresponding to each scale choice are
labeled as \textsc{sd}, \textsc{mad}, \textsc{$M$-scale}. For each
scale, we consider four estimators, the raw estimators where no
smoothing is used, the estimators obtained by penalizing the scale
function defined in (\ref{estFPCsmooth1}), those obtained by penalizing
the norm defined in (\ref{estFPCsmooth2}) and the sieve estimators
defined in (\ref{sieveFPC}). In all tables, as in\vspace*{1pt} Section
\ref{prop}, the $j$th principal direction estimators related to each
method are labeled as $\wphi_{\raw,j}$, $\wphi_{\smooths,j}$, $\wphi
_{\smoothn,j}$ and $\wphi_{\sieve,j}$, respectively.

When using the penalized estimators, several values for the penalizing
parameters $\tau$ and $\rho$ were chosen. Since large values of the
smoothing parameters make the penalizing term the dominant component
regardless of the amount of contamination considered, we choose $\tau$
and $\rho$ equal to $a n^{-\alpha}$ for $\alpha=3$ and $4$ and $a$
equal to $0.05$, $0.10$, $0.25$, $0.5$, $0.75$, $1$, $1.5$ and $2$.

For the sieve estimators based on the Fourier basis, ordered as
$\{1, \cos( \pi x)$, $\sin( \pi x),\ldots, \cos( q_n \pi x), \sin(
q_n \pi x), \ldots\}$, the values $p_n=2 q_n$ with $q_n=5$, $10$ and
$15$ are used, while for the sieve estimators based on the $B$-splines,
the dimension of the linear space considered is selected as $p_n= 10$,
$20$ and $50$. The basis for the $B$-splines is generated from the R
function $\mathit{cSplineDes}$, with the knots being equally spaced in the
interval $[-1,1]$ and
the number of knots equal to $p_n+1$.
To conserve space, we only report here the results corresponding to
$p_n=30$ and $p_n=50$ for the Fourier and $B$-spline basis,
respectively. More extensive simulation results are listed in the
technical report [\citet{bali2}].

%s8.3 #&#
\subsection{Simulation settings}\label{setting}

The sample was generated using a finite Karhu\-nen--Lo\`eve expansion
with the functions, $\phi_i\dvtx[-1,1] \rightarrow\mathbb{R}$,
$i=1,2,3$, where
$ \phi_1(x)= \sin( 4 \pi x)$, $\phi_2(x) = \cos( 7 \pi x)$ and
$\phi
_3(x)= \cos( 15 \pi x)$.
It is worth noticing that, when considering the sieve estimators based
on the Fourier basis, the third component cannot be detected when $q_n
< 15$, since in this case~$\phi_3(x)$ is orthogonal to the estimating
space. Likewise, the second component cannot be detected when $q_n < 7$.

We performed $NR=1\mbox{,}000$ replications generating independent samples $\{
X_i\}_{i=1}^n$ of size $n= 100$ following the model $X_i = Z_{i1} \phi
_1 + Z_{i2} \phi_2 + Z_{i3} \phi_3$, where $\bZ
_{i}=(Z_{i1},Z_{i2},Z_{i3})\trasp$ are independent vectors whose
distribution will depend on the situation to be considered. The central
model, denoted $C_0$, corresponds to Gaussian samples. We also consider
four contaminations of the central model, labeled $C_{2}$, $C_{3,a}$,
$C_{3,b}$ and $C_{23}$ depending on the components to be contaminated.
Contamination models are commonly considered in robust statistics since
they tend be the more difficult models to be robust against and are the
basis for the concept of bias robustness, see \citet{maro} for
further discussion. In all these situations $Z_{i1}, Z_{i2}$ and
$Z_{i3}$ are also independent. For each of the models, we took $\sigma
_1=4$, $\sigma_2=2$ and $\sigma_3=1$. The central model and the
contaminations can be described as follows:\vspace*{8pt}

${C_{0}}$: $Z_{i1} \sim N( 0,\sigma_1^2)$, $Z_{i2} \sim N(0,\sigma
_2^2)$ and $Z_{i3} \sim N(0,\sigma_3^2)$.

${C_{2}}$: $Z_{i2}$ are independent and identically distributed as $
0.8 N(0,\sigma_2^2)+\break 0.2 N(10,\allowbreak 0.01)$, while
$Z_{i1} \sim N(0,\sigma_1^2)$ and $Z_{i3} \sim N(0,\sigma_3^2)$. This
contamination corresponds to a strong contamination on the second
component and changes the mean value of the generated data $Z_{i2}$ and
also the first principal component. Note that $\var(Z_{i2})=19.202$.
%%queda phi2, phi1, phi3

${C_{3,a}}$: $Z_{i1} \sim N( 0,\sigma_1^2)$, $Z_{i2} \sim N(0,\sigma
_2^2)$ and $Z_{i3} \sim0.8 N(0,\sigma_3^2)+ 0.2
N(15$, $0.01)$. This contamination corresponds to a strong
contamination on the third component. Note that $\var(Z_{i3})=36.802$.
%%%queda phi3, phi1, phi2

${C_{3,b}}$: $Z_{i1} \sim N( 0,\sigma_1^2)$, $Z_{i2} \sim N(0,\sigma
_2^2)$ and $Z_{i3} \sim0.8 N(0,\sigma_3^2)+ 0.2
N(6, 0.01)$. This contamination corresponds to a strong
contamination on the third component. Note that $\var(Z_{i3})=6.562$.
%%%queda phi1, phi3, phi2

${C_{23}}$: $Z_{ij}$ are independent and such that $Z_{i1} \sim
N(0,\sigma_1^2)$, $Z_{i2} \sim0.9 N(0,\sigma_2^2) + 0.1
N(15, 0.01)$ and $Z_{i3} \sim0.9 N(0,\sigma_3^2) +
0.1 N(20, 0.01)$. This contamination corresponds to a mild
contamination on the last two components. Note that $\var
(Z_{i2})=23.851$ and $\var(Z_{i3})=36.901$.\vspace*{8pt} %queda phi3, phi2, phi1

We also include a long-tailed model, namely a Cauchy model,
labe-\break led~${C_{\mathrm{Cauchy}}}$,~which is defined by taking $\bZ_i\,{\sim}\,\mathcal
{C}_3(0,\bSi)$ with $\bSi\,{=}\,\operatorname{diag}(\sigma_1^2,\sigma
_2^2,\sigma
_3^2)$, where $\mathcal{C}_p(0,\bSi)$ stands for the $p$-dimensional
elliptical Cauchy distribution centered at $0$ with scatter matrix
$\bSi
$. For this situation, the covariance operator does not exist, and thus
the classical principal components are not defined.

It is worth noting that the directions $\phi_1$, $\phi_2$ and $\phi_3$
correspond to the classical principal components for the case $C_0$,
but not
necessarily for the other cases. For instance, when $\sigma_{\rob
}^2=\operatorname{VAR}$, $C_{3,a}$ interchanges the order between~$\phi_1$
and $\phi_3$, that is, $\phi_3=\phi_{\rob,1}(C_{3,a})$ as defined in
(\ref{MAXROB}), and so it corresponds to the first principal component
of the covariance operator, while $\phi_1$ is the second and $\phi_2$
is the third one.

%s8.4 #&#
\subsection{Simulation results}\label{resultados}
For each situation, we compute the estimators of the first three
principal directions and the square distance between the true and the
estimated direction (normalized to have $L^2$ norm 1), that is,
\[
D_j=\biggl\|\frac{\wphi_j}{\| \wphi_j\|}-\phi_j\biggr\|^2 .
\]
Note that all the estimators except those penalizing the norm, are such
that \mbox{$\| \wphi_j\|=1$}. Tables 4 to
9 of the supplementary
material [\citet{bali3b}] report the means of $D_j$ over replications,
which hereafter is referred to as mean square error. To help understand
the influence of the grid size $m$ on the estimators, Tables
3, 4 and 5 give the mean squared errors
for $m=50, 100, 150, 250$ and $250$, under $C_0$ for various values of
the penalizing parameters. As can be seen, for the first two components
some slight improvement is observed when using $m=250$ as opposed to
$m=50$ points, but at the expense of increasing the computing time
about 2.6 fold. On the other hand, for the third principal direction,
taking $m=100$ compared to $m=50$ reduces the mean square error by at
least a half for the penalized estimators,
while the gain is not so prominent for the raw estimates. The size
$m=50$ is selected for presentation in the remainder of our study since it
provides a reasonable a compromise between the performance of the
estimators and the computational time.

To help understand the effect of penalization, consider Table
6. This table shows results for the raw and
penalized estimators under $C_0$ for different choices of the
penalizing parameters. From this table, we see that a better
performance is achieved in most cases when $\alpha=3$ is used. To be
more precise, the results in Table~6 show
that the best choice for $\wphi_{\smooths,j}$ is $\rho=2 n^{-3}$ for
all $j$. Note that $\rho=1.5 n^{-3}$ gives fairly similar results when
using the $M$-scale, reducing the mean squared error relative to the
raw estimate by about a half and a third for $j=2$ and $3$, respectively.

For the norm penalized estimators, $\wphi_{\smoothn,j}$, the best
choice for the penalizing parameter seems to depend upon both the
component to be estimated and the scale estimator used. For instance,
when using the standard deviation, the best choice is $\tau=0.10
n^{-3}$, for $j=1$ and $2$ while for $j=3$ a smaller order is needed to
obtain an improvement over the raw estimators, with the value $\tau
=0.75 n^{-4}$ leading to a small gain over the raw estimators. For the
robust procedures, larger values are needed to see an advantage to
using the penalized norm approach relative to the raw estimators. For
example, for $j=1$, the largest reduction is observed when $\tau=2
n^{-3}$ while for $j=2$, the best choices correspond to $\tau=0.5
n^{-3}$ and $\tau= 0.25 n^{-3}$ when using the \textsc{mad} and
$M$-scale, respectively. When using the $M$-scale, a good compromise is
to choose $\tau=0.75 n^{-3} $, which gives a reduction of around 30\%
and 50\% for the first and second principal directions, respectively,
although smaller values of~$\tau$ are again better for estimating the
third component.

Based upon\vspace*{1pt} the above observations, we report here only the results
corresponding to $\rho=1.5 n^{-3}$ and $\tau=0.75 n^{-3}$ for the
penalized estimators $\wphi_{\smoothn,j}$ and $\wphi_{\smooths,j}$,
respectively, under the contamination models and the Cauchy model.
Results for other choices of $\rho$ and $\tau$ are given in
\citet{bali2}.

The simulation study confirms the expected inadequate behavior of the
classical estimators in the presence of outliers. Under contamination,
the classical estimators of the principal directions do not estimate
the target directions very accurately. This is also the case when
considering the Cauchy distribution. Curiously, though, the principal
directions, under the Cauchy model, do not seem to be totally arbitrary
and they partially recover $\phi_1$, $\phi_2$ and $\phi_3$ when the
standard deviation is used, although not as well as when using a robust
scale, even though the covariance operator does not exist, nor do the
population principal directions as defined in~(\ref{MAX}).

The robust estimators of the first principal directions are not heavily
affected by any of the contaminations, while the estimates of the
second and third principal directions appear to be most affected under
model $C_{3,a}$. In particular, for the third direction, the
projection-pursuit estimators based on an the \textsc{mad} seems to be
most affected by this type of contamination when penalizing the norm,
although much less so than the classical methods. With respect to the
contamination model~$C_{3,a}$, the estimators $\wphi_{\smoothn,j}$,
which are the robust penalized norm estimators, tend to have the best
performance among all the robust competitors for the first two
components, and, in particular, when using the $M$-scale; see Table
8. It is worth noting that the
classical estimators of the first component are not affected by this
contamination when penalizing the norm since the penalization dominates
the contaminated variances. The same phenomena is observed under
$C_{3,b}$ when using the classical estimators for the selected amount
of penalization. For the raw estimators, the sensitivity of the
classical estimators under this contamination can be observed in Table
7. We refer to \citet{bali2}
for the behavior when other values of the smoothing parameters
are chosen.

As noted in \citet{s}, for the classical estimators, some
degree of smoothing in the procedure based on penalizing the norm will
give a better estimation of $\phi_j$ in the $L^2$ sense under mild
conditions. In particular, both the procedure penalizing the norm and
the scale provide some improvement with respect to the raw estimators
if $\Psi(\phi_j)<\Psi(\phi_\ell)$, when $j<\ell$. This means that the
principal directions are rougher as the eigenvalues decrease [see
\citet{pez} and \citet{s}], which is also
reflected in our simulation study. The advantages of the smooth
projection pursuit procedures are most striking when estimating $\phi
_2$ and $\phi_3$ with an $M$-scale and using the penalized scale approach.

As expected, when using the sieve estimators, the Fourier basis gives
the best performance over all the methods under $C_0$, since our data
were generated from this basis; see Table 9. The choice of the $B$-spline
basis give similar results quite to those obtained with $\wphi
_{\smooths
,j}$ when estimating the first direction, except under $C_{\mathrm{Cauchy}}$
where the penalized estimators show a~better performance. For the
second and third components, the estimators obtained with the
$B$-spline basis show larger the mean square errors than the raw or
penalized estimators.

%s8.5 #&#
\subsection{$K$th fold simulation}\label{kfold}
Table \ref{tabtiempos} reports the computing times in minutes for 1,000
replications and for a fixed value of $\tau$ or $\rho$, run on a
computer Core Quad I7~930 (2.80 GHz) with 8 Gb of Ram memory. We also
report the computing times when using the sieve approach with the
Fourier basis and a~fixed value of $p_n$. This suggests that the
leave-one-out cross-validation may be difficult to perform, and so a
$K$-fold approach is adopted instead. It is worth noticing that the
robust procedures based on the \textsc{mad} are much faster than those
based on the $M$-scale, so they may be preferred in terms of computing
time. However, as mentioned in Section~\ref{prelim}, the main
disadvantage of the \textsc{mad} is its low efficiency and lack of
smoothness, which is related to the discontinuity of its influence
function. This is particularly important when estimating principal
components in the finite-dimensional setting, since as it was pointed
out by \citet{cui} and \citet{crouruiz2} the
variances of some elements of the estimated principal directions may
blow up when using the \textsc{mad} leading to highly inefficient
estimators. As expected and mentioned in Section \ref{resultados},
Table 6 reveals a high loss of efficiency
for the \textsc{mad}, in our functional setting, for any choice of the
smoothing parameter.

%
%t1 #&#
%
\begin{table}
\tablewidth=250pt
\caption{Computing times in minutes for 1,000
replications and a fixed value of $\tau$, $\rho$ or $p_n$ (when using
the Fourier basis)}\label{tabtiempos}
\begin{tabular*}{\tablewidth}{@{\extracolsep{\fill}}lccc@{}}
\hline
& \textbf{\textsc{sd}} & \textbf{\textsc{mad}} & \textbf{\textsc{$\bolds{M}$-scale}}\\
\hline
$\wphi_{\raw}$ & \hphantom{0}5.62 & \hphantom{0}6.98 & 17.56\\[2pt]
$\wphi_{\smooths}$ & \hphantom{0}7.75 & \hphantom{0}9.00 & 20.18\\[2pt]
$\wphi_{\smoothn}$ & 31.87 & 33.21 & 44.04\\[2pt]
$\wphi_{\sieve}$ & \hphantom{0}0.5\hphantom{0} & \hphantom{0}5.22 & 16.08 \\
\hline
\end{tabular*}
\end{table}

For the procedure which penalizes the scale or the norm, the smoothing
parameters $\rho$ and $\tau$ are selected using the procedure described
in Section \ref{smoothpar} with $K=4$ and $\ell=1$. Due to the
extensive computing time, we have only performed $500$ replications.
The simulation results when penalizing the scale function, that is, for
the estimators defined through (\ref{estFPCsmooth1}), are reported in
Table \ref{tabtablameanL2squarenormKfold}. Under $C_0$, when estimating
the second and third principal directions, the robust estimators based
on the $M$-scale combined with a penalization in the scale clearly have
smaller mean square error than the raw estimators, while those
penalizing the norm improve the performance of the raw estimators and
also that of $\wphi_{\smooths,j}$, on the first and second
directions.\looseness=1

%
%t2 #&#
%
\begin{table}
\caption{Mean values of $\|\wphi_j/\|\wphi_j\|-\phi_j\|^2$
when the penalizing parameter is selected using
$K$-fold cross-validation}\label{tabtablameanL2squarenormKfold}
\begin{tabular*}{\tablewidth}{@{\extracolsep{\fill}}lccccccc@{}}
\hline
& & \multicolumn{3}{c}{$\bolds{\wphi_{\smooths,j}}$} & \multicolumn
{3}{c}{$\bolds{\wphi_{\smoothn,j}}$} \\[-4pt]
& & \multicolumn{3}{c}{\hrulefill} & \multicolumn
{3}{c@{}}{\hrulefill}\\
\textbf{Model} & \textbf{Scale estimator}
& $\bolds{j=1}$ & $\bolds{j=2}$ & $\bolds{j=3}$ & $\bolds{j=1}$
& $\bolds{j=2}$ & $\bolds{j=3}$ \\
\hline
$C_0$ & \textsc{SD} & 0.0073 & 0.0094 & 0.0078 & 0.0075 & 0.0094 & 0.0360 \\
& \textsc{mad} & 0.0662 & 0.0993 & 0.0634 & 0.0497 & 0.0660 &
0.2573 \\
& $M$-scale & 0.0225 & 0.0311 & 0.0172 & 0.0208 & 0.0271 & 0.0839 \\
[4pt]
$C_2$ & \textsc{SD} & 1.2840 & 1.2837 & 0.0043 & 1.2076 & 1.2232 & 0.0301 \\
& \textsc{mad} & 0.3731 & 0.3915 & 0.0504 & 0.3360 & 0.3770 &
0.2832 \\
& $M$-scale & 0.4261 & 0.4286 & 0.0153 & 0.3679 & 0.4049 & 0.1607 \\
[4pt]
$C_{3,a}$ & \textsc{SD} & 1.7840 & 1.8901 & 1.9122 & 1.7795 & 1.8861 & 1.9134 \\
& \textsc{mad} & 0.2271 & 0.5227 & 0.5450 & 0.0573 & 0.2289 &
0.9540 \\
& $M$-scale & 0.2176 & 0.4873 & 0.5437 & 0.0257 & 0.1187 & 0.8710 \\
[4pt]
$C_{3,b}$ & \textsc{SD} & 0.0192 & 0.8350 & 0.8525 & 0.0173 & 0.5902 & 0.7502 \\
& \textsc{mad} & 0.0986 & 0.3930 & 0.3820 & 0.0553 & 0.1417 &
0.5167 \\
& $M$-scale & 0.0404 & 0.2251 & 0.2285 & 0.0241 & 0.1080 & 0.3174 \\
[4pt]
$C_{23}$ & \textsc{SD} & 1.7645 & 0.5438 & 1.6380 & 1.7537 & 0.6496 & 1.4305 \\
& \textsc{mad} & 0.2407 & 0.3443 & 0.2064 & 0.1414 & 0.2214 &
0.6824 \\
& $M$-scale & 0.2613 & 0.3707 & 0.2174 & 0.1313 & 0.1870 & 0.5901 \\
[4pt]
$C_{\mathrm{Cauchy}}$ & \textsc{SD} & 0.3580 & 0.4835 & 0.2287 & 0.2862 & 0.3525 & 0.3435 \\
& \textsc{mad} & 0.0788 & 0.1511 & 0.1082 & 0.0613 &
0.0855 & 0.3147 \\
& $M$-scale & 0.0444 & 0.0707 & 0.0434 & 0.0349 & 0.0463 & 0.1465 \\
\hline
\end{tabular*}
\end{table}

From the results in Table \ref{tabtablameanL2squarenormKfold} we
observe that the classical estimators are sensitive to the
contaminations in the simulation settings, and, except for
contaminations in the third component, the robust counterpart shows a~clear advantage. Note that $C_{3,b}$ affects more the classical
estimators when the smoothing parameter is selected by the robust
$K$-fold cross-validation method than for the fixed values studied in
the previous section. This can be explained by the fact that
contamination $C_{3,b}$ is a mild contamination in the third component
which has a large $\|\phi^{\prime\prime}_3\|^2$, and so the classical
estimators are more sensitive to it, just as the raw estimators, if
smaller values of the smoothing parameter are chosen. It is worth
noticing that the penalized robust estimators based on the $M$-scale
improve the performance of the raw estimators based on the $M$-scale,
even under contaminations, when the penalizing parameter is selected
using the $K$-fold approach. This advantage is more striking when
penalizing the norm and when the two first principal components are considered.

Note that we choose $\ell=1$, and so our focus was on the first
principal component. To improve the observed performance, a different
approach should be considered, maybe by selecting a different smoothing
parameter for each principal direction.

%s9 #&#
\section{Concluding remarks}\label{concl}
In this paper, we propose robust principal component analysis for
functional data based on a projection-pursuit approach. The different
procedures correspond to robust versions of the unsmoothed principal
component estimators, to the estimators obtained penalizing the scale
and to those obtained by penalizing the norm. A sieve approach based on
approximating the elements of the unit ball by elements over
finite-dimensional spaces is also considered. In particular, the
procedures based on smoothing and sieves are new. A robust
cross-validation procedure is introduced to select the smoothing parameters.
Consistency results are derived for the four type of estimators.
Finally, a simulation study confirms the expected inadequate behavior
of the classical
estimators in the presence of outliers, with the robust procedures
performing significantly better. In particular, the procedure\vadjust{\goodbreak} based on
an $M$-scale combined with a penalization in the norm, where the
smoothing parameter is selected via a robust $K$-fold cross-validation,
is recommended.

Among other contributions we highlight the following:\vspace*{8pt}

(a) We obtain the continuity of the principal directions and
eigenvalue functionals, which implies the asymptotic qualitative
robustness of the corresponding estimators. This extends the results of
\citet{li} from Euclidean spaces to infinite-dimensional
Hilbert spaces, where the unit ball is not compact. Noncompactness
poses technical challenges which we overcome with tools from functional
analysis.

(b) Our results not only include the finite-dimensional case but also
improve upon some of the results obtained in that situation for the
projection pursuit estimators. For example, the assumptions in
\citet{li} regarding the robust scale functional are stronger than
ours. Also, to derive the consistency of the raw estimators, we only
require uniform convergence over the unit ball of $s_n(\alpha)$
to~$\sigma(\alpha)$, which holds if the scale functional
$\sigma_{\rob}$ is continuous. This improves upon the consistency
results given in \citet{cui}, who require a uniform Bahadur
expansion for $s_n(\alpha )$ over the unit ball.

(c) A key step in proving the continuity of the projection pursuit
functional is to show that weak convergence of probability measures
over a~Hilbert space implies uniform convergence of the laws of the
projections of the stochastic processes, that is, Theorem
\ref{theorem62}. This uniform convergence result can be useful in other
statistical problems where projection methods are considered.

(d) The proofs for the penalized estimators include, as particular
cases, the estimators defined by \citet{ris} and studied
by \citet{pez}, and those considered by \citet{s}. Extending the results to scale estimators other than the
standard deviation required more challenging arguments since, unlike
the classical setting, the projection-pursuit index $s^2_n(\alpha)$
cannot be expressed in the simple form $\langle\alpha, \bGa_n \alpha
\rangle$ for some compact operator $\bGa_n$.

\begin{appendix}\label{appenA}
%s10 #&#
\section*{Appendix}
In this Appendix, we give the proofs of the results stated in Section
\ref{consist}. Some technicalities are omitted, and we refer to the
technical report [\citet{bali2}] for details.
Before presenting the proof, some additional notation is needed.

Denote by $\mathcal{L}_{m-1}$ the linear space spanned by $\{\phi
_{\rob
,1},\ldots,\phi_{\rob, m-1}\}$, and let~$\widehat{\mathcal
{L}}_{m-1}$ be
the linear space spanned by the first $m-1$ estimated principal
directions, that is, by $\{\wphi_{\sieve,1},\ldots,\wphi_{\sieve
,m-1}\}$
or $\{\wphi_{1},\ldots,\wphi_{m-1}\}$, where it will be clear in each
case which linear space we are considering. The latter includes the
situation of the linear spaces spanned by $\{\wphi_{\raw,1},\ldots
,\wphi
_{\raw, m-1}\}$, $\{\wphi_{\smooths,1},\ldots,\wphi_{\smooths
,m-1}\}$
and $\{\wphi_{\smoothn,1},\ldots,\wphi_{\smoothn,m-1}\}$. Finally, for
any linear space~$\mathcal L$, $\pi_\mathcal{L}\dvtx\mathcal{H}\to
\mathcal
{L}$ stands for the orthogonal projection onto the linear space
$\mathcal L$, which exists if $\mathcal L$ is a closed linear space. In
particular, $\pi_{\mathcal{L}_{m-1}}$, $\pi_{\widehat{\mathcal
{L}}_{m-1}}$ and~$\pi_{\mathcal{H}_{p_n}}$ are well defined.

Moreover, for the sake of simplicity, denote by $\mathcal{T}_k =
\mathcal{L}_{k}^{\bot}$ the linear space orthogonal to $\phi_1,\ldots
,\phi_k$ and by $\pi_k = \pi_{\mathcal{T}_k}$ the orthogonal projection
with respect to the inner product defined in $\mathcal H$. On the other
hand, let $\widehat{\pi}_{\tau,k}$ be the projection onto the linear
space orthogonal to $\wphi_{1},\ldots, \wphi_{k}$ in the space
$\mathcal
{H}_{\smooth}$ in the inner product $\langle\cdot, \cdot\rangle
_{\tau
}$, that is, for any $\alpha\in\mathcal{H}_{\smooth}$,
$\widehat{\pi}_{\tau,k} (\alpha) = \alpha- \sum_{j=1}^k \langle
\alpha,
\wphi_{ j} \rangle_{\tau}\wphi_{ j}$.
Moreover, let $\widehat{\mathcal{T}}_{\tau,k}$ stand for the linear
space orthogonal to $\widehat{\mathcal{L}}_k$ with the inner product
$\langle\cdot, \cdot\rangle_{\tau}$. Thus, $ \widehat{\pi}_{\tau,k}$
is the orthogonal projection onto $\widehat{\mathcal{T}}_{\tau,k}$ with
respect to this inner product.

\begin{pf*}{Proof of Theorem \ref{theorem41}}
First note that the
fact that $\sigma_{\rob}$ is a scale functional entails that $\sigma
_n(\alpha)=\|\alpha\| \sigma_n(\alpha/\|\alpha\| )$. Thus from
assumption (ii) and the fact that $\|\alpha\| \leq\|\alpha\|_\tau$, we
get that
%
%e10.1 #&#
%
\begin{equation}\label{S1enmenorque1}\quad
\sup_{\|\alpha\|\le1}|\sigma_n^2(\alpha) - \sigma^2(\alpha)
|\convpp0 \quad\mbox{and}\quad \sup_{\|\alpha\|_{\tau}\le1}
|\sigma_n^2(\alpha) - \sigma^2(\alpha)|\convpp0 .
\end{equation}

(a) To prove that $\wlam_1 \convpp\sigma^2(\phi_{\rob,1})$ it is
enough to show that
%
%e10.2 #&#
%e10.3 #&#
%
\begin{eqnarray}
\label{SMdesi1}
\sigma^2(\phi_{\rob,1}) &\geq& \wlam_1 + o_{\mathrm{a.s.}}(1) ,
\\
\label{SMdesi2}
\sigma^2(\phi_{\rob,1}) &\leq& \wlam_1 + o_{\mathrm{a.s.}}(1),
\end{eqnarray}
where $ o_{\mathrm{a.s.}}(1)$ stands for a term converging to $0$
almost surely.

Note that from (\ref{S1enmenorque1}), we get that $a_{n,1}=\sigma
_n^2(\wphi_1)-\sigma^2(\wphi_1)\convpp0$ and $b_{n,1}=\sigma
_n^2(\phi
_{\rob,1})-\sigma^2(\phi_{\rob,1})\convpp0$.
Using that $\sigma$ is a scale functional and that\break $\sigma^2(\phi
_{\rob
,1}) = \sup_{\alpha\in\mathcal{S}} \sigma^2(\alpha)$, we obtain
easily that
\[
\sigma^2(\phi_{\rob,1}) \geq\sigma^2 \biggl( \frac{\wphi_1}{\|\wphi
_1\|
}\biggr)=\frac{\sigma^2 (\wphi_1)}{\|\wphi_1\|^2}\geq\sigma^2(\wphi_1)=
\sigma_n^2(\wphi_1) - a_{n,1}= \wlam_1 + o_{\mathrm{a.s.}}(1)
\]
concluding the proof of (\ref{SMdesi1}).

To derive (\ref{SMdesi2}), note that since $\phi_{\rob,1}\in
\mathcal
{H}_{\smooth}$, $\| \phi_{\rob,1} \|_{\tau}<\infty$ and $\|\phi
_{\rob
,1}\|_{\tau} \geq\|\phi_{\rob,1}\|=1$,
then, defining\vspace*{1pt} $\beta_1=\phi_{\rob,1}/ \|\phi_{\rob,1}\|_{\tau}$, we
have that $\|\beta_1\|_{\tau}=1$, which implies that $
\wlam_1 = \sigma_n^2(\wphi_1) \geq\sigma_n^2(\wphi_1)-\rho\Psi
(\wphi
_1)\ge\sigma_n^2(\beta_1)-\rho\Psi(\beta_1)$.
Hence, using that $\sigma_{\rob}$ is a scale functional and that
$\Psi
(a\alpha)=a^2\Psi(\alpha)$, for any $a\in\real$, we get
\begin{eqnarray*}
\wlam_1 &\ge& \sigma_n^2(\beta_1)-\rho\Psi(\beta_1)= \frac
{\sigma
_n^2(\phi_{\rob,1})-\rho\Psi(\phi_{\rob,1})}{\|\phi_{\rob,1}\|
_{\tau}^2}
= \frac{\sigma^2(\phi_{\rob,1})+b_{n,1}-\rho\Psi(\phi_{\rob
,1})}{\|\phi
_{\rob,1}\|_{\tau}^2}.
\end{eqnarray*}
When $\rho=0$, we have defined $\rho\Psi(\phi_{\rob,1})= 0$ and
similarly when $\tau=0$. So from now on, we will assume that $\tau_n>0$
and $\rho_n>0$.
Since\vadjust{\goodbreak} $b_{n,1}=o_{\mathrm{a.s.}}(1)$, $\rho\convpp0$ and $\tau
\convpp0$, we have that $\rho\Psi(\phi_{\rob,1})\convpp0$ and $\|
\phi
_{\rob,1}\|_{\tau} \convpp\|\phi_{\rob,1}\|= 1$, concluding the proof
of (\ref{SMdesi2}). Hence, $\wlam_1 \convpp\sigma^2(\phi_{\rob,1})$.

From (\ref{S1enmenorque1}) and the fact that $\|\wphi_1\| \leq1$, we
obtain that $\wlam_1-\sigma^2(\wphi_1)=\sigma_n^2(\wphi_1)-\sigma
^2(\wphi_1)\convpp0$. Therefore, using that $\wlam_1 \convpp\sigma
^2(\phi_{\rob,1})$, we get that
%
%e10.4 #&#
%
\begin{equation}\label{sigmawphi}
\sigma^2(\wphi_1) \convpp\sigma^2(\phi_{\rob,1}) .
\end{equation}
Moreover, the inequalities $\sigma^2(\phi_{\rob,1}) \geq\sigma^2 (
{\wphi_1}/{\|\wphi_1\|})\geq\sigma^2(\wphi_1)$
obtained above also imply that
%
%e10.5 #&#
%
\begin{equation}\label{sigmawphi2}
\sigma^2 ( {\wphi_1}/{\|\wphi_1\|}) \convpp\sigma^2(\phi_{\rob ,1}) .
\end{equation}
Using that $\|\wphi_1\|_{\tau}=1$, we get that $\tau\Psi(\wphi_1)
=1-\|\wphi_1\|^2=1-{\sigma^2(\wphi_1)}/\sigma^2(\wphi_1 /\allowbreak\|\wphi_1\|)$.
Hence, (\ref{sigmawphi}) and (\ref{sigmawphi2}) entail that
$\tau
\Psi(\wphi_1)\convpp0$.

It only remains to show that $\rho\Psi(\wphi_1)\convpp0$, which
follows easily from the fact that $\wlam_1\convpp\sigma^2(\phi
_{\rob
,1})$, $\sigma_n^2(\phi_{\rob,1})\convpp\sigma^2(\phi_{\rob
,1})$, $\rho
\convpp0$\vspace*{1pt} and $\|\phi_{\rob,1}\|_{\tau} \convpp1$ since
$\wlam_1 \geq\sigma_n^2(\wphi_1) - \rho\lceil\wphi_1, \wphi_1
\rceil\geq
(\sigma_n^2(\phi_{\rob,1}) - \rho\lceil\phi_{\rob,1}, \phi
_{\rob,1}
\rceil)/\|\phi_{\rob,1}\|_{\tau}^2$.

Note that we have not used the weak continuity of $\sigma$ as a
function of $\alpha$ to derive~(a).

(b) Note that since $\|\wphi_{1}\|_{\tau}=1$, we have that $\|
\wphi
_{1}\|\le1$. Moreover, from~(a), $\|\wphi_{1}\| \convpp1$. Let
$\wtphi
_1=\wphi_1/\|\wphi_1\|$, then $\wtphi_1\in\mathcal{S}$ and $\sigma
(\wtphi_1)=\sigma(\wphi_1)/\|\wphi_1\|$. Using that $\sigma
^2(\wphi
_{1})\convpp\sigma^2(\phi_{\rob,1})$ and $\|\wphi_{1}\| \convpp
1$, we
obtain that $\sigma^2(\wtphi_{1})\convpp\sigma^2(\phi_{\rob,1})$, and
thus the proof follows using Lemma \ref{lemma41}.

(c) Let us show that $\wlam_m \convpp\sigma^2(\phi_{\rob,m})$.
The proof will be done in several steps by showing
%
%e10.6 #&#
%e10.7 #&#
%e10.8 #&#
%
\begin{eqnarray}
\label{lema}
&\displaystyle {\sup_{\|\alpha\|_{\tau} \leq1} }|\sigma^2(\pi_{m-1} \alpha) -
\sigma
_n^2(\widehat{\pi}_{\tau,m-1} \alpha)|
\convpp0,&
\\
\label{desilamm1}
&\sigma^2(\phi_{\rob,m}) \geq\wlam_m + o_{\mathrm{a.s.}}(1),&
\\
\label{desigenm}
&\sigma^2(\phi_{\rob,m}) \leq\wlam_m + o_{\mathrm{a.s.}}(1) .&
\end{eqnarray}
Note that (\ref{lema}) corresponds to an extension of assumption (ii)
while (\ref{desilamm1}) and (\ref{desigenm}) are analogous to (\ref
{SMdesi1}) and (\ref{SMdesi2}).

We begin by proving (\ref{lema}). Note that $ {\sup_{\|\alpha\|_{\tau}
\leq1}} |\sigma_n^2(\pi_{m-1} \alpha) -\sigma
^2(\widehat{\pi}_{\tau,m-1} \alpha)| \leq{\sup_{\|\alpha\|_{\tau}
\leq1}} |\sigma^2(\pi_{m-1} \alpha) - \sigma^2(\widehat{\pi}_{\tau,m-1}
\alpha)| +\vspace*{1pt} {\sup_{\|\alpha\|_{\tau} \leq1}}
|\sigma_n^2(\widehat{\pi}_{\tau,m-1} \alpha) -
\sigma^2(\widehat{\pi}_{\tau,m-1} \alpha)|$. Using
(\ref{S1enmenorque1}) and the fact that if $\|\alpha\|_{\tau} \leq1$,
then $\|\widehat{\pi}_{\tau,m-1} \alpha\|_{\tau} \leq1$, we get that
the second term on the right-hand side converges to 0 almost surely. To
complete the proof of (\ref{lema}), it remains to show that
%
%e10.9 #&#
%
\begin{equation}\label{falta}
{\sup_{\|\alpha\|_{\tau} \leq1}} |\sigma^2(\pi_{m-1} \alpha) -
\sigma
^2(\widehat{\pi}_{\tau,m-1} \alpha)|\convpp0 .
\end{equation}
Using that $\wphi_j \convpp\phi_{\rob,j}$ and that $\tau\Psi(
\wphi
_j)=\tau\lceil\wphi_j,\wphi_j \rceil\convpp0$, for $1\le j\le m-1$
and arguing as in \citet{s} [see \citet{bali2} for
details], we get that,\vadjust{\goodbreak} for $1\le j \le m-1$,
$
{\sup_{\|\alpha\|_{\tau} \leq1}} \|\langle\alpha, \phi_{\rob
,j}\rangle
\phi_{\rob,j}-\langle\alpha, \wphi_j\rangle_{\tau}\wphi_j\|
\convpp0
$,
entailing that ${\sup_{\|\alpha\|_{\tau} \leq1}} \|\widehat{\pi
}_{\tau
,m-1} \alpha-\pi_{m-1} \alpha\|\convpp0$.
Therefore, using that $\sigma$ is weakly uniformly continuous over the
unit ball, we get easily that (\ref{falta}) holds, concluding the proof
of (\ref{lema}).

As in (a), we will next show that (\ref{desilamm1}) holds. Using again
that $\sigma$ is a scale functional, we get easily that
$\sup_{\alpha\in\mathcal{S}\cap\mathcal{T}_{m-1}} \sigma^2(\alpha) =
\sup_{\alpha\in\mathcal{S}} \sigma^2(\pi_{m-1} \alpha) $, so using
again that $\|\wphi_m\|\le\|\wphi_m\|_\tau=1$, we obtain\vspace*{1pt} that
$\sigma^2(\phi_{\rob,m}) =\break\sup_{\alpha\in\mathcal{S}} \sigma^2(\pi
_{m-1} \alpha) \geq\sigma^2(\pi_{m-1}
{\wphi_m}/{\|\wphi_m\|})\ge\sigma^2 (\pi_{m-1} {\wphi_m})$. From
(\ref{lema}) and the fact that $\|\wphi_m\|_{\tau}=1$, we get that
$b_m=\sigma^2(\pi_{m-1} \wphi_m)-\sigma_n^2 (\widehat{\pi}_{\tau,m-1}
\wphi_m)\convpp0$, and so since $\widehat{\pi}_{\tau,m-1}
{\wphi_m}=\wphi_m$ and $\|\wphi_m\|\le1$, we get that $
\sigma^2(\phi_{\rob,m}) \geq\sigma^2(\pi_{m-1} \wphi_m)= \sigma_n^2
(\widehat{\pi}_{\tau,m-1} \wphi_m) + o_{\mathrm{a.s.}}(1) = \wlam_m +
o_{\mathrm{a.s.}}(1)$, completing the proof of (\ref{desilamm1}).

We will show now (\ref{desigenm}). Note that
$\phi_{\rob,m}\in\mathcal{H}_{\smooth}$, so that $\| \phi_{\rob,m}
\|_{\tau}<\infty$ and $\| \phi_{\rob,m}\|_{\tau}\to\| \phi_{\rob,m}
\|=1$. Using\vspace*{2pt} that $\sigma_{\rob}$ is a scale functional, the fact that
$\wlam_m = \sigma_n^2(\wphi_m)\ge\sigma_n^2(\wphi_m) -\rho\Psi(\wphi
_m)=\sup_{\|\alpha\|_{\tau}=1, \alpha\in\widehat{\mathcal{T}}_{\tau
,m-1}}\{ \sigma_n^2(\alpha) -\rho\Psi(\alpha)\}$ and that for any
$\alpha\in\mathcal{H}_{\smooth}$ such that $\|\alpha\|_{\tau}=1$ we
have that $\| \widehat{\pi}_{\tau,m-1} \alpha\|_{\tau}\le1$, we get
easily that $\wlam_m \ge
\sup_{\|\alpha\|_{\tau}=1}\{\sigma_n^2(\widehat{\pi}_{\tau,m-1} \alpha)
-\rho\Psi(\widehat{\pi}_{\tau,m-1} \alpha)\}$, and so $\wlam_m
\geq({\sigma_n^2( \widehat{\pi}_{\tau,m-1}\phi_{\rob,m} )- \rho\Psi(
\widehat{\pi}_{\tau,m-1} \phi_{\rob,m})}
)/{\|\phi_{\rob,m}\|_{\tau}^2}$. From (\ref{lema}) we obtain that
$d_m=\sigma_n^2(\widehat{\pi}_{\tau,m-1}
\phi_{\rob,m})-\sigma^2(\pi_{m-1} \phi_{\rob,m})\convpp0$. Moreover,
the fact that $\tau\convpp0$ entails that
$\|\phi_{\rob,m}\|_{\tau}\convpp\|\phi_{\rob,m}\|=1$. On the other
hand, using that $\rho\Psi(\wphi_{\ell})\convpp0$, $1\le\ell\le m-1$,
and the fact that $\rho\convpp0$ implies that $\rho\Psi(\phi_{\rob,m})
= o_{\mathrm{a.s.}}(1)$, analogous arguments to those considered in
\citet{pez} allow us to show that $
\rho\Psi(\wpi_{m-1}\phi_{\rob,m})=\break\rho\lceil\wpi_{m-1}\phi_{\rob,m},
\wpi_{m-1}\phi_{\rob,m} \rceil\convpp0$. Hence, we get that
\begin{eqnarray*}
\wlam_m %&\geq& \sigma_n^2(\wphi_m) -\rho\Psi(\wphi_m) \ge\frac{
&\geq&\frac{\sigma^2( {\pi}_{ m-1}\phi_{\rob,m}
)+d_m- \rho
\Psi( \widehat{\pi}_{\tau,m-1} \phi_{\rob,m})}{1+o(1)}
\\
&\geq&\frac{\sigma^2( {\pi}_{ m-1}\phi_{\rob,m})+d_m-
o_{\mathrm{a.s.}}(1)}{1+o(1)}\\
&=&\sigma^2( \phi_{\rob,m}
)+o_{\mathrm{a.s.}}(1) ,
\end{eqnarray*}
where the last equality follows from the fact that ${\pi}_{ m-1}\phi
_{\rob,m}=\phi_{\rob,m}$.

Therefore, from (\ref{desilamm1}) and (\ref{desigenm}), we obtain that
$\wlam_m \convpp\sigma^2(\phi_{\rob,m} )$.

On the other hand, (\ref{lema}) entails that $\wlam_m- \sigma
^2(\wphi
_m)=\sigma_n^2(\wphi_m) - \sigma^2(\wphi_m)\convpp0$, which together
with $\wlam_m\convpp\sigma^2(\phi_{\rob,m})$ implies that $ \sigma
^2(\wphi_m)\convpp\sigma^2(\phi_{\rob,m})$.

To complete the proof of (c), it remains to show that $\tau\Psi(\wphi
_m) \convpp0$ and $\rho\Psi(\wphi_m) \convpp0$. As in (a), we have
that the following inequalities converge to equalities:
%
%e10.10 #&#
%
\begin{equation}
\label{cotasigmaphim}\qquad
\sigma^2(\phi_{\rob,m}) \geq\sigma^2\biggl(\pi_{m-1} \frac{\wphi
_m}{\|
\wphi_m\|}\biggr)\geq\sigma^2(\pi_{m-1} \wphi_m)=
\wlam_m +
o_{\mathrm{a.s.}}(1) .\vadjust{\goodbreak}
\end{equation}
Using\vspace*{1pt} that $\sigma$ is a scale estimator and that
$\|\wphi_m\|_\tau=1$, we get that $\tau\Psi(\wphi_m) = 1 -
\|\wphi_m\|^2 = 1 -{\sigma^2(\pi _{m-1}\wphi
_m)}/{\sigma^2(\pi_{m-1}\wphi_m /\|\wphi_m\|)} $, which together with
(\ref{cotasigmaphim}) entails that the second term on the right-hand
side is $1+ o_{\mathrm{a.s.}}(1)$ and so, $\tau \lceil\wphi_m,\wphi_m
\rceil\convpp0$, entailing that $\|\wphi _m\| \convpp1$.

On the other hand, we also have that
%
%e10.11 #&#
%
\begin{equation}
\label{cotasigmawphim}
\wlam_m=\sigma_n^2(\wphi_{m}) \geq\sigma_n^2(\wphi_m
)-\rho
\Psi(\wphi_m)\geq\sigma^2(\phi_{\rob,m}) + o_{\mathrm{a.s.}}(1) ,
\end{equation}
so using that $\wlam_m=\sigma_n^2(\wphi_{m})\convpp\sigma^2
(\phi
_{\rob,m})$, we obtain that $\rho\Psi(\wphi_m)\convpp0$,
concluding the proof of (c).

(d) We have already proved that when $m=1$ the result holds. We
proceed by induction and assume\vspace*{1pt} that $\langle\wphi_\ell, \phi_{\rob
,\ell} \rangle^2 \to1$, $\tau\Psi(\wphi_{ \ell}) \convpp0$ and
$\rho
\Psi(\wphi_{ \ell}) \convpp0$ for $1\le\ell\le m-1$, to show that
$\langle\wphi_m, \phi_{\rob,m} \rangle^2 \to1$. Without loss of
generality, we can assume that $\wphi_\ell\convpp\phi_{\rob,\ell}$,
for $1\le\ell\le m-1$. Denote by $\wtphi_j={\wphi_j}/{\|\wphi_j\|}$.
Then, for $1 \le\ell\le m-1$, $\|\wphi_\ell\|\to1$, and so $\wtphi
_\ell\convpp\phi_{\rob,\ell}$. It suffices to show that $\langle
\phi
_{\rob,m}, \wtphi_m\rangle^2\convpp1$.

Using (c) we have that $\sigma^2(\wphi_m) \convpp\sigma^2(\phi
_{\rob,
m})$ and that $\|\wphi_m\|\convpp1$, and so $\sigma^2(\wtphi_m)
\convpp\sigma^2(\phi_{\rob, m})$. The proof follows now from Lemma
\ref{lemma41} if we show that $\langle\wtphi_m, \wtphi_\ell\rangle
\convpp
0$, $1\le\ell\le m-1$.

Using that $\tau\Psi(\wphi_\ell) \convpp0$, for $1\le\ell\le m-1$,
and that from (c) $\tau\Psi(\wphi_m) \convpp0$ we get that $\tau
\lceil\wphi_\ell,\wphi_m \rceil\convpp0$ for $1\le\ell\le m-1$.
Therefore, the fact that $\langle\wphi_m, \wphi_\ell\rangle_\tau=0$
entails that
$\langle\wphi_m, \wphi_\ell\rangle=\langle\wphi_m, \wphi_\ell
\rangle
_\tau- \tau\lceil\wphi_\ell,\wphi_m \rceil\convpp0$, and so
$\langle
\wtphi_m, \wtphi_\ell\rangle\convpp0$, concluding the proof.
\end{pf*}
\end{appendix}

\section*{Acknowledgments}

We wish to thank the Associate Editor and three anonymous referees for
valuable comments which led to an improved version of the original
paper.

\begin{supplement}[id=suppA]
\sname{Supplement A}
\stitle{Robust functional principal components\\}
\slink[doi]{10.1214/11-AOS923SUPPA} %[doi,text={...}] - jei reikia
%suskaldyti doi
\sdatatype{.pdf}
\sfilename{aos923\_suppa.pdf}
\sdescription{In this Supplement, we give the proof of some of the
results stated in Sections \ref{consist} and \ref{appen}.}
\end{supplement}

\begin{supplement}[id=suppB]
\sname{Supplement B}
\stitle{Robust functional principal components\\}
\slink[doi]{10.1214/11-AOS923SUPPB} %[doi,text={...}] - jei reikia
%suskaldyti doi
\sdatatype{.pdf}
\sfilename{aos923\_suppb.pdf}
\sdescription{In this Supplement, we report the results obtained
in the Monte Carlo study for the raw estimators and for the penalized
ones when the smoothing parameters are fixed.}
\end{supplement}

% imsref loaded by lrinkeviciute, 2011-11-30 10:16:27
% imsref loaded by lrinkeviciute, 2011-11-30 10:42:49
%

\printaddresses

\end{document}